\documentclass[leqno,draft]{article}

\newtheorem{theorem}{Theorem}
\newtheorem{lemma}[theorem]{Lemma}
\newtheorem{proposition}[theorem]{Proposition}
\newtheorem{definition}[theorem]{Definition}
\newtheorem{corollary}[theorem]{Corollary}

\newcommand{\begintheorem}{\addtocounter{equation}{1}\begin{theorem}}
\newcommand{\beginlemma}{\addtocounter{equation}{1}\begin{lemma}}
\newcommand{\beginproposition}{\addtocounter{equation}{1}\begin{proposition}}
\newcommand{\begindefinition}{\addtocounter{equation}{1}\begin{definition}}
\newcommand{\begincorollary}{\addtocounter{equation}{1}\begin{corollary}}

\begin{document}

\title{Some remarks about solenoids}

\author{Stephen Semmes \\
        Rice University}

\date{}

\maketitle

\begin{abstract}
A basic family of solenoids is discussed, especially from the point of
view of analysis on metric spaces.
\end{abstract}

\tableofcontents

\part{Basic examples}
\label{basic examples}

\section{A Cartesian product}
\label{cartesian product}
\setcounter{equation}{0}

        Let $r \ge 2$ be an integer, and consider the Cartesian product
\begin{equation}
\label{X = prod_{l = 0}^infty ({bf R} / r^l {bf Z})}
        X = \prod_{l = 0}^\infty ({\bf R} / r^l \, {\bf Z}).
\end{equation}
More precisely, the real line ${\bf R}$ is a commutative group with
respect to addition, $r^l \, {\bf Z}$ is the discrete subgroup of
${\bf R}$ consisting of integer multiples of $r^l$, and ${\bf R} / r^l
\, {\bf Z}$ is the corresponding quotient group.  The quotient ${\bf
R} / r^l \, {\bf Z}$ may also be considered as a compact Hausdorff
topological space and a $1$-dimensional smooth manifold in the usual
way.  The quotients ${\bf R} / r^l {\bf Z}$ may actually be considered
as Lie groups, because the group operations are given by smooth
mappings.  These Lie groups are all isomorphic to each other, and to
the multiplicative group of complex numbers with modulus equal to $1$.
The Cartesian product $X$ consists of the sequences $x = \{x_l\}_{l =
0}^\infty$ with $x_l \in {\bf R} / r^l \, {\bf Z}$, and is a compact
Hausdorff space with respect to the product topology.  Of course, $X$
is a commutative group as well, where the group operations are defined
coordinatewise.  It is easy to see that the group operations on $X$
are continuous with respect to the product topology on $X$, so that
$X$ is a topological group.

        Because $r^{l + 1} \, {\bf Z} \subseteq r^l \, {\bf Z}$, there
is a natural homomorphism from ${\bf R} / r^{l + 1} \, {\bf Z}$ onto
${\bf R} / r^l \, {\bf Z}$ for each $l \ge 0$.  An element $x =
\{x_l\}_{l = 0}^\infty$ of $X$ is said to be a \emph{coherent
sequence} if $x_l$ is the image in ${\bf R} / r^l \, {\bf Z}$ of $x_{l
+ 1} \in {\bf R} / r^{l + 1} \, {\bf Z}$ for each $l$.  Note that the
set $Y$ of coherent sequences in $X$ is a closed subgroup of $X$ with
respect to the topology and group structure described in the previous
paragraph.

        Let $q_l$ be the usual quotient mapping from ${\bf R}$ onto
${\bf R} / r^l \, {\bf Z}$ for each $l$.  Consider the mapping $q$
from ${\bf R}$ into $X$ defined by
\begin{equation}
\label{q(a) = {q_l(a)}_{l = 0}^infty}
        q(a) = \{q_l(a)\}_{l = 0}^\infty
\end{equation}
for each $a \in {\bf R}$.  This is a continuous homomorphism from
${\bf R}$ into $X$ with trivial kernel.  Observe that $q(a)$ is a
coherent sequence in $X$ for each $a \in {\bf R}$, because $q_l$ is
the same as the composition of $q_{l + 1}$ with the natural mapping
from ${\bf R} / r^{l + 1} \, {\bf Z}$ onto ${\bf R} / r^l \, {\bf Z}$
for each $l$.  Thus $q({\bf R}) \subseteq Y$, and in fact $q({\bf R})$
is dense in $Y$, so that $Y$ is the same as the closure of $q({\bf
R})$ in $X$.  To see this, let $x \in Y$ and a positive integer $L$ be
given, and choose $a \in {\bf R}$ such that $q_L(a) = x_L$.  The
coherence condition implies that $q_l(a) = x_l$ for each $l \le L$,
and hence that $q(a)$ is arbitrarily close to $x$ with respect to the
product topology on $X$, as desired.  It follows that $Y$ is
connected, since it is the closure of the connected set $q({\bf R})$.

        Let $\pi_l$ be the $l$th coordinate projection of $X$ onto
${\bf R} / r^l \, {\bf Z}$, so that
\begin{equation}
\label{pi_l(x) = x_l}
        \pi_l(x) = x_l
\end{equation}
for each $x \in X$ and $l \ge 0$.  Thus $\pi_l$ is a continuous
homomorphism from $X$ onto ${\bf R} / r^l \, {\bf Z}$, and similarly
the restriction of $\pi_l$ to $Y$ defines a continuous homomorphism
from $Y$ onto ${\bf R} / r^l \, {\bf Z}$ for each $l$.  By
construction, the restriction of $\pi_l$ to $Y$ is the same as the
composition of the restriction of $\pi_{l + 1}$ to $Y$ with the
natural homomorphism from ${\bf R} / r^{l + 1} \, {\bf Z}$ onto ${\bf
R} / r^l \, {\bf Z}$.  If $y \in Y$ is in the kernel of $\pi_0$, then
it follows that
\begin{equation}
\label{pi_l(y) in {bf Z} / r^l {bf Z}}
        \pi_l(y) \in {\bf Z} / r^l \, {\bf Z}
\end{equation}
for each $l \ge 0$.  In particular, the kernel of the restriction of
$\pi_0$ to $Y$ is totally disconnected.

\section{A nice metric}
\label{nice metric}
\setcounter{equation}{0}

        Let $\phi_l$ be the standard isomorphism between ${\bf R} /
r^l \, {\bf Z}$ and the unit circle ${\bf T}$ in the complex plane
${\bf C}$.  Thus
\begin{equation}
\label{phi_l(q_l(a)) = exp (2 pi i r^{-l} a)}
        \phi_l(q_l(a)) = \exp (2 \, \pi \, i \, r^{-l} \, a)
\end{equation}
for every $a \in {\bf R}$, where $\exp z$ is the usual complex
exponential function on ${\bf C}$.  Remember that
\begin{equation}
\label{|exp (i t)| = 1}
        |\exp (i \, t)| = 1
\end{equation}
for every $t \in {\bf R}$, where $|\zeta|$ denotes the modulus of
$\zeta \in {\bf C}$.  Note that
\begin{equation}
\label{d_l(x_l, y_l) = |phi_l(x_l) - phi_l(y_l)|}
        d_l(x_l, y_l) = |\phi_l(x_l) - \phi_l(y_l)|
\end{equation}
defines a metric on ${\bf R} / r^l \, {\bf Z}$, and that the topology
on ${\bf R} / r^l \, {\bf Z}$ determined by this metric is the same as
the quotient topology corresponding to the standard topology on ${\bf
R}$.  This is the same as saying that $\phi_l$ is a homeomorphism from
${\bf R} / r^l \, {\bf Z}$ onto ${\bf T}$ with respect to the topology
on ${\bf T}$ induced by the standard Euclidean metric on ${\bf C}$.

        If $x, y \in X$, then put
\begin{equation}
\label{d(x, y) = max_{l ge 0} r^{-l} |phi_l(x_l) - phi_l(y_l)|}
        d(x, y) = \max_{l \ge 0} r^{-l} \, |\phi_l(x_l) - \phi_l(y_l)|.
\end{equation}
Of course,
\begin{equation}
\label{|phi_l(x_l) - phi_l(y_l)| le |phi_l(x_l)| + |phi_l(y_l)| = 2}
        |\phi_l(x_l) - \phi_l(y_l)| \le |\phi_l(x_l)| + |\phi_l(y_l)| = 2,
\end{equation}
which implies that
\begin{equation}
\label{lim_{l to infty} r^{-l} |phi_l(x_l) - phi_l(y_l)| = 0}
        \lim_{l \to \infty} r^{-l} \, |\phi_l(x_l) - \phi_l(y_l)| = 0
\end{equation}
for every $x, y \in X$.  This ensures that the maximum in (\ref{d(x,
y) = max_{l ge 0} r^{-l} |phi_l(x_l) - phi_l(y_l)|}) is always attained.

        It is easy to see that $d(x, y)$ satisfies the requirements of
a metric on $X$.  In particular, the triangle inequality for $d(x, y)$
can be verified using the triangle inequality for (\ref{d_l(x_l, y_l)
= |phi_l(x_l) - phi_l(y_l)|}) for each $l$.  The topology on $X$
corresponding to $d(x, y)$ is the same as the product topology
discussed in the previous section.  More precisely,
\begin{equation}
\label{d(x, y) < t}
        d(x, y) < t
\end{equation}
for some positive real number $t$ if and only if
\begin{equation}
\label{r^{-l} |phi_l(x_l) - phi_l(y_l)| < t}
        r^{-l} \, |\phi_l(x_l) - \phi_l(y_l)| < t
\end{equation}
for each $l \ge 0$ such that $2 \, r^{-l} \ge t$.  Thus (\ref{d(x, y)
< t}) only involves finitely many coordinates of $x$ and $y$ for any
given $t > 0$, which implies that open subsets of $X$ with respect to
$d(x, y)$ are also open with respect to the product topology.
Conversely, one can show that open subsets of $X$ with respect to the
product topology also also open with respect to $d(x, y)$, because
(\ref{d(x, y) < t}) implies that any finite number of coordinates of
$x$ and $y$ are arbitrarily close to each other when $t$ is
sufficiently small.  Of course, we are using the fact that
(\ref{d_l(x_l, y_l) = |phi_l(x_l) - phi_l(y_l)|}) determines the
quotient topology on ${\bf R} / r^l \, {\bf Z}$ corresponding to the
standard topology on ${\bf R}$ for each $l \ge 0$ here.

        If $x_l, y_l, z_l \in {\bf R} / r^l \, {\bf Z}$, then
\begin{eqnarray}
\label{d_l(x_l + z_l, y_l + z_l) = d_l(x_l, y_l)}
 d_l(x_l + z_l, y_l + z_l) & = & |\phi_l(x_l + z_l) - \phi_l(y_l + z_l)| \\
 & = & |\phi_l(x_l) \, \phi_l(z_l) - \phi_l(y_l) \, \phi_l(z_l)| \nonumber \\
 & = & |\phi_l(x_l) - \phi_l(y_l)| \, |\phi_l(z_l)| \nonumber \\
 & = & |\phi_l(x_l) - \phi_l(y_l)| = d_l(x_l, y_l). \nonumber
\end{eqnarray}
This shows that $d_l(x_l, y_l)$ is invariant under translations on
${\bf R} / r^l \, {\bf Z}$ for each $l \ge 0$.  It follows that
\begin{equation}
\label{d(x + z, y + z) = d(x, y)}
        d(x + z, y + z) = d(x, y)
\end{equation}
for every $x, y, z \in X$, so that $d(x, y)$ is also invariant under
translations on $X$.

\section{Another Cartesian product}
\label{another cartesian product}
\setcounter{equation}{0}

        Consider the Cartesian product
\begin{equation}
\label{X_0 = prod_{l = 1}^infty ({bf Z} / r^l {bf Z})}
        X_0 = \prod_{l = 1}^\infty ({\bf Z} / r^l \, {\bf Z}).
\end{equation}
Thus the elements of $X_0$ are sequences $x = \{x_l\}_{l = 1}^\infty$
such that $x_l \in {\bf Z} / r^l \, {\bf Z}$ for each $l$.  We can
identify $X_0$ with a subset of $X$, because ${\bf Z} / r^l \, {\bf Z}
\subseteq {\bf R} / r^l \, {\bf Z}$ for each $l \ge 1$, and by
extending $x = \{x_l\}_{l = 1}^\infty$ to $l = 0$ by taking $x_0 = 0$
in ${\bf R} / {\bf Z}$.  Note that $X_0$ corresponds to a closed
subgroup of $X$ with respect to coordinatewise addition in this way.
The topology on $X_0$ induced by the product topology on $X$ is the
same as the product topology on $X_0$ that corresponds to taking the
discrete topology on ${\bf Z} / r^l \, {\bf Z}$ for each $l$.
Actually, $r^l \, {\bf Z}$ is an ideal in the ring of integers for
each $l$, so that each quotient ${\bf Z} / r^l \, {\bf Z}$ may be
considered as a commutative ring.  It follows that $X_0$ is a
commutative ring with respect to coordinatewise addition and
multiplication as well.  It is easy to see that multiplication on
$X_0$ is continuous with respect to the product topology, so that
$X_0$ is a topological ring.

        As before, there is a natural ring homoorphism from ${\bf Z} /
r^{l + 1} \, {\bf Z}$ onto ${\bf Z} / r^l \, {\bf Z}$ for each $l \ge
1$, because $r^{l + 1} \, {\bf Z} \subseteq r^l \, {\bf Z}$.  An
element $x = \{x_l\}_{l = 1}^\infty$ of $X_0$ is said to be a
\emph{coherent sequence} if $x_l$ is the image in ${\bf Z} / r^l \,
{\bf Z}$ of $x_{l + 1} \in {\bf Z} / r^{l + 1} \, {\bf Z}$ for each
$l$.  Thus $x$ is a coherent sequence in $X_0$ if and only if the
corresponding element of $X$ is a coherent sequence in the sense of
Section \ref{cartesian product}.  Equivalently, the set $Y_0$ of
coherent sequences in $X_0$ can be identified with the subset of $X$
which is the intersection of the set $Y$ of coherent sequences in $X$
with the subset of $X$ identified with $X_0$.  Note that $Y_0$ is a
closed subring of $X_0$.

        Let $\widetilde{q}_l$ be the natural quotient mapping from
${\bf Z}$ onto ${\bf Z} / r^l \, {\bf Z}$ for each $l \ge 1$.  This is
the same as the restriction of the quotient mapping $q_l : {\bf R} \to
{\bf R} / r^l \, {\bf Z}$ from Section \ref{cartesian product} to
${\bf Z}$, although now $\widetilde{q}_l$ is a ring homomorphism from
${\bf Z}$ onto ${\bf Z} / r^l \, {\bf Z}$.  Similarly, let
$\widetilde{q}$ be the mapping from ${\bf Z}$ into $X_0$ defined by
\begin{equation}
\label{widetilde{q}(a) = {widetilde{q}_l(a)}_{l = 1}^infty}
        \widetilde{q}(a) = \{\widetilde{q}_l(a)\}_{l = 1}^\infty
\end{equation}
for each $a \in {\bf Z}$.  This is a ring homomorphism from ${\bf Z}$
into $X_0$ with trivial kernel, and which is the same as the
restriction of the embedding $q : {\bf R} \to X$ defined in Section
\ref{cartesian product} to ${\bf Z}$ when we identify $X_0$ with a
subset of $X$ as before.  In particular, $\widetilde{q}(a)$ is a
coherent sequence in $X_0$ for each $a \in {\bf Z}$, for the same
reasons as before.  One can also check that $\widetilde{q}({\bf Z})$
is dense in $Y_0$, so that $Y_0$ is the same as the closure of
$\widetilde{q}({\bf Z})$ in $X_0$ with respect to the product
topology.  Of course, $X_0$ is obviously totally disconnected, and so
$Y_0$ is too.

        Let $\pi_0$ be the $l = 0$ coordinate projection of $X$ onto
${\bf R} / {\bf Z}$, as in Section \ref{cartesian product}.  The
kernel of the restriction of $\pi_0$ to $Y$ consists of the coherent
sequences $y = \{y_l\}_{l = 0}^\infty$ in $X$ such that $y_0 = 0$ in
${\bf R} / {\bf Z}$.  Because of the coherence condition, this implies
that $y_l \in {\bf Z} / r^l \, {\bf Z}$ for each $l \ge 1$.  Thus the
kernel of the restriction of $\pi_0$ to $Y$ corresponds exactly to the
subset of $X$ identified with $Y_0$.

\section{Another metric}
\label{another metric}
\setcounter{equation}{0}

        Let $x$ and $y$ be distinct elements of the set $X_0$ defined
in the previous section, and let $l(x, y)$ be the smallest positive
integer $l$ such that $x_l \ne y_l$.  Equivalently, $l(x, y)$ is the
largest positive integer $l$ such that $x_j = y_j$ for every $j < l$.
Put
\begin{equation}
\label{rho(x, y) = r^{-l(x, y) + 1}}
        \rho(x, y) = r^{-l(x, y) + 1}.
\end{equation}
If $x = y$, then we put $\rho(x, y) = 0$, which corresponds to taking
$l(x, y) = +\infty$ in (\ref{rho(x, y) = r^{-l(x, y) + 1}}).  Of
course,
\begin{equation}
\label{l(x, y) = l(y, x)}
        l(x, y) = l(y, x)
\end{equation}
for every $x, y \in X_0$, which implies that
\begin{equation}
\label{rho(x, y) = rho(y, x)}
        \rho(x, y) = \rho(y, x).
\end{equation}
Similarly,
\begin{equation}
\label{l(x, z) ge min(l(x, y), l(y, z))}
        l(x, z) \ge \min(l(x, y), l(y, z))
\end{equation}
for every $x, y, z \in X_0$, and hence
\begin{equation}
\label{rho(x, z) le max(rho(x, y), rho(y, z))}
        \rho(x, z) \le \max(\rho(x, y), \rho(y, z)).
\end{equation}
It follows that $\rho(x, y)$ defines an ultrametric on $X_0$, which
means that $\rho(x, y)$ is a metric on $X_0$ that satisfies the
stronger ultrametric version (\ref{rho(x, z) le max(rho(x, y), rho(y,
z))}) of the triangle inequality.  

        It is easy to see that the topology on $X_0$ determined by
$\rho(x, y)$ is the same as the product topology corresponding to the
discrete topology on each factor ${\bf R} / r^l \, {\bf Z}$ in
(\ref{X_0 = prod_{l = 1}^infty ({bf Z} / r^l {bf Z})}).  We also have
that
\begin{equation}
\label{l(x + z, y + z) = l(x, y)}
        l(x + z, y + z) = l(x, y)
\end{equation}
for every $x, y, z \in X_0$, so that
\begin{equation}
\label{rho(x + z, y + z) = rho(x, y)}
        \rho(x + z, y + z) = \rho(x, y).
\end{equation}
Thus $\rho(x, y)$ is invariant under translations on $X_0$.

        We would like to compare this metric with the one in Section
\ref{nice metric}.  As before, $x, y \in X_0$ may be identified with
elements of $X$, by taking $x_0 = y_0 = 0$ in ${\bf R} / {\bf Z}$.  In
this case, (\ref{d(x, y) = max_{l ge 0} r^{-l} |phi_l(x_l) -
phi_l(y_l)|}) reduces to
\begin{equation}
\label{d(x, y) = max_{l ge 1} r^{-l} |phi_l(x_l) - phi_l(y_l)|}
        d(x, y) = \max_{l \ge 1} r^{-l} \, |\phi_l(x_l) - \phi_l(y_l)|.
\end{equation}
We may as well suppose that $x \ne y$, since otherwise $d(x, y) =
\rho(x, y) = 0$, so that (\ref{d(x, y) = max_{l ge 1} r^{-l}
|phi_l(x_l) - phi_l(y_l)|}) reduces further to
\begin{equation}
\label{d(x, y) = max_{l ge l(x, y)} r^{-l} |phi_l(x_l) - phi_l(y_l)|}
        d(x, y) = \max_{l \ge l(x, y)} r^{-l} \, |\phi_l(x_l) - \phi_l(y_l)|.
\end{equation}
In particular,
\begin{equation}
\label{d(x, y) le 2 r^{-l(x, y)} = 2 r^{-1} rho(x, y)}
        d(x, y) \le 2 \, r^{-l(x, y)} = 2 \, r^{-1} \, \rho(x, y),
\end{equation}
by (\ref{|phi_l(x_l) - phi_l(y_l)| le |phi_l(x_l)| + |phi_l(y_l)| =
2}).

        In the other direction, we can take $l = l(x, y)$ in (\ref{d(x, y) = 
max_{l ge l(x, y)} r^{-l} |phi_l(x_l) - phi_l(y_l)|}), to get that
\begin{equation}
\label{d(x, y) ge r^{-l(x, y)} ...}
        d(x, y) \ge r^{-l(x, y)} \,
                  |\phi_{l(x, y)}(x_{l(x, y)}) - \phi_{l(x, y)}(y_{l(x, y)})|.
\end{equation}
Under these conditions, $x_{l(x, y)}$ and $y_{l(x, y)}$ are distinct
elements of ${\bf Z} / r^{l(x, y)} \, {\bf Z}$, and hence
\begin{equation}
\label{|phi_{l(x, y)}(x_{l(x, y)}) - phi_{l(x, y)}(y_{l(x, y)})| ge ..., X_0}
        |\phi_{l(x, y)}(x_{l(x, y)}) - \phi_{l(x, y)}(y_{l(x, y)})|
                               \ge |\exp (2 \, \pi \, i \, r^{-l(x, y)}) - 1|.
\end{equation}
If $x, y \in Y_0$, so that $x$ and $y$ are coherent sequences, then
$x_{l(x, y)}$ and $y_{l(x, y)}$ are distinct elements of ${\bf Z} /
r^{l(x, y)} \, {\bf Z}$ which are equal module $r^{l(x, y) - 1} \,
{\bf Z}$, and
\begin{equation}
\label{|phi_{l(x, y)}(x_{l(x, y)}) - phi_{l(x, y)}(y_{l(x, y)})| ge ..., Y_0}
        |\phi_{l(x, y)}(x_{l(x, y)}) - \phi_{l(x, y)}(y_{l(x, y)})|
                               \ge |\exp (2 \, \pi \, i \, r^{-1}) - 1|.
\end{equation}
Combining this with (\ref{d(x, y) ge r^{-l(x, y)} ...}), we get that
\begin{equation}
\label{d(x, y) ge r^{-1} |exp (2 pi i r^{-1}) - 1| rho(x, y)}
 d(x, y) \ge r^{-1} \, |\exp (2 \, \pi \, i \, r^{-1}) - 1| \, \rho(x, y)
\end{equation}
for every $x, y \in Y_0$.

\section{$r$-Adic integers}
\label{r-adic integers}
\setcounter{equation}{0}

        Let $a$ be a nonzero integer, and let $l(a)$ be the largest
nonnegative integer $l$ such that $a$ is an integer multiple of $r^l$.
If $b$ is another nonzero integer, then it is easy to see that
\begin{equation}
\label{l(a + b) ge min(l(a), l(b))}
        l(a + b) \ge \min(l(a), l(b))
\end{equation}
and
\begin{equation}
\label{l(a b) ge l(a) + l(b)}
        l(a \, b) \ge l(a) + l(b).
\end{equation}
The $r$-adic absolute value $|a|_r$ of $a$ is defined by
\begin{equation}
\label{|a|_r = r^{-l(a)}}
        |a|_r = r^{-l(a)}.
\end{equation}
Of course, we put $|a|_r = 0$ when $a = 0$, which corresponds to
taking $l(a) = +\infty$.  Thus we get that
\begin{equation}
\label{|a + b|_r le max(|a|_r, |b|_r)}
        |a + b|_r \le \max(|a|_r, |b|_r)
\end{equation}
and
\begin{equation}
\label{|a b|_r le |a|_r |b|_r}
        |a \, b|_r \le |a|_r \, |b|_r
\end{equation}
for all integers $a$, $b$.  The $r$-adic metric on ${\bf Z}$ is defined by
\begin{equation}
\label{delta_r(a, b) = |a - b|_r}
        \delta_r(a, b) = |a - b|_r.
\end{equation}
It is easy to see that this defines a metric on ${\bf Z}$, and more
precisely an ultrametric on ${\bf Z}$, since
\begin{equation}
\label{delta_r(a, c) le max(delta_r(a, b), delta_r(b, c))}
        \delta_r(a, c) \le \max(\delta_r(a, b), \delta_r(b, c))
\end{equation}
for every $a, b, c \in {\bf Z}$, by (\ref{|a + b|_r le max(|a|_r, |b|_r)}).

        Let $a$, $b$ be integers, and let $\widetilde{q}(a)$,
$\widetilde{q}(b)$ be their images in $X_0$, as in Section
\ref{another cartesian product}.  We would like to check that
\begin{equation}
\label{rho(widetilde{q}(a), widetilde{q}(b)) = delta_r(a, b)}
        \rho(\widetilde{q}(a), \widetilde{q}(b)) = \delta_r(a, b),
\end{equation}
where $\rho(x, y)$ is the ultrametric on $X_0$ defined in Section
\ref{another metric}.  To do this, it suffices to show that
\begin{equation}
\label{l(widetilde{q}(a), widetilde{q}(b)) - 1 = l(a - b)}
        l(\widetilde{q}(a), \widetilde{q}(b)) - 1 = l(a - b),
\end{equation}
where $l(x, y)$ is defined for $x, y \in X_0$ as in the previous
section.  Thus $l(\widetilde{q}(a), \widetilde{q}(b))$ is the smallest
positive integer $l$ such that $\widetilde{q}_l(a) \ne \widetilde{q}_l(b)$,
which is the same as saying that $l(\widetilde{q}(a), \widetilde{q}(b)) - 1$
is the largest nonnegative integer $k$ such that $\widetilde{q}_j(a) = 
\widetilde{q}_j(b)$ for every $j \le k$.  Remember that $\widetilde{q}_j$
is the natural quotient homomorphism from ${\bf Z}$ onto ${\bf Z} / r^j \, 
{\bf Z}$, so that $\widetilde{q}_j(a) = \widetilde{q}_j(b)$ exactly when
$a - b$ is an integer multiple of $r^j$.  It follows that $l(\widetilde{q}(a),
\widetilde{q}(b)) - 1$ is the same as the largest nonnegative integer $k$
such that $a - b$ is an integer multiple of $r^k$, which is also the same
as $l(a - b)$, as desired.  Note that we could have reduced to the case
where $b = 0$ at the beginning of the argument, because $\widetilde{q}$
is a homomorphism from ${\bf Z}$ into $X_0$, and because of the 
translation-invariance of the metrics involved.

        A sequence $x(1) = \{x_l(1)\}_{l = 1}^\infty, x(2) =
\{x_l(2)\}_{l = 1}^\infty, x(3) = \{x_l(3)\}_{l = 1}^\infty, \ldots$
of elements of $X_0$ converges to an element $x = \{x_l\}_{l =
1}^\infty$ of $X_0$ with respect to the product topology discussed in
Section \ref{another cartesian product}, or equivalently with respect
to the ultrametric $\rho(\cdot, \cdot)$, if and only if for each
positive integer $n$ we have that $x_l(n) = x_l$ for all sufficiently
large $l$, depending on $n$.  Similarly, if $x(1), x(2), x(3), \ldots$
is a Cauchy sequence in $X_0$ with respect to $\rho(\cdot, \cdot)$,
then it is easy to see that $x_l(n)$ is eventually constant in $l$ for
each $n$, and hence that $x(1), x(2), x(3), \ldots$ converges in
$X_0$.  This shows that $X_0$ is complete as a metric space with
respect to $\rho(\cdot, \cdot)$, which could also be derived from the
compactness of $X_0$.  It follows that $Y_0$ is complete as a metric
space with respect to $\rho(\cdot, \cdot)$ too, because $Y_0$ is a
closed subset of $X_0$.

        Thus $Y_0$ can be identified with the completion ${\bf Z}_r$
of ${\bf Z}$ with respect to the $r$-adic metric, since
$\widetilde{q}$ is an isometric embedding of ${\bf Z}$ onto a dense
subset of $Y_0$, and $Y_0$ is complete with respect to $\rho(\cdot,
\cdot)$.  In particular, the ring structure on $Y_0$ defined by
coordinatewise addition and multiplication corresponds to the ring
structure on ${\bf Z}_r$ obtained by extending addition and
multiplication on ${\bf Z}$ to ${\bf Z}_r$ by continuity.  The
completion ${\bf Z}_r$ of ${\bf Z}$ with respect to the $r$-adic
metric is known as the ring of $r$-adic integers, especially when $r =
p$ is a prime number.  In this case, equality holds in (\ref{l(a b) ge
l(a) + l(b)}) and (\ref{|a b|_r le |a|_r |b|_r}), and the $p$-adic
absolute value and metric can be defined on the field ${\bf Q}$ of
rational numbers.  The completion ${\bf Q}_p$ of ${\bf Q}$ with
respect to the $p$-adic metric is known as the field of $p$-adic
numbers, and ${\bf Z}_p$ is the same as the closure of ${\bf Z}$ in
${\bf Q}_p$.

\section{A nice mapping}
\label{nice mapping}
\setcounter{equation}{0}

        Consider the mapping $A$ from ${\bf R} \times Y_0$ into $Y$ defined by
\begin{equation}
\label{A(a, x) = q(a) + x}
        A(a, x) = q(a) + x.
\end{equation}
Remember that $q$ maps ${\bf R}$ into $Y$ as in Section \ref{cartesian
product}, and that we identify $x = \{x_l\}_{l = 1}^\infty \in Y_0$
with an element of $Y$ by setting $x_0 = 0$ in ${\bf R} / {\bf Z}$.
Thus (\ref{A(a, x) = q(a) + x}) is defined by taking the sum of $q(a)$
and $x$ as elements of $Y$ as a subgroup of $X$ as a commutative group
with respect to coordinatewise addition.  More precisely, $q$ is a
homomorphism of ${\bf R}$ into $Y$ with respect to addition, and hence
$A$ is a homomorphism from ${\bf R} \times Y_0$ into $Y$ with respect
to coordinatewise addition on ${\bf R} \times Y_0$.

        Suppose that $(a, x) \in {\bf R} \times Y_0$ is in the kernel
of $A$, so that $q(a) + x = 0$ in $Y$.  In particular, the $l = 0$
coordinate of $q(a) + x$ is equal to $0$ in ${\bf R} / {\bf Z}$, which
implies that $q_0(a) = 0$ in ${\bf R} / {\bf Z}$, because $x \in Y_0$.
It follows that $a \in {\bf Z}$, and that $x = \widetilde{q}(-a)$ in
$Y_0$ in the notation of Section \ref{another cartesian product}.
Conversely, if $a \in {\bf Z}$ and $x = \widetilde{q}(-a)$ in $Y_0$,
then $A(a, x) = 0$.

        Let $y = \{y_l\}_{l = 0}^\infty$ be any element of $Y$.  If
$y_0 = 0$, then $y$ can be identified with an element of $Y_0$, and
$y$ is in the image of $A$.  Otherwise, we can choose $a \in {\bf R}$
such that $q_0(a) = y_0$ in ${\bf R} / {\bf Z}$, so that the $l = 0$
coordinate of $y - q(a) \in Y$ is equal to $0$.  This implies that $y
- q(a)$ corresponds to an element of $Y_0$, and hence that $y = q(a) +
(y - q(a))$ is in the image of $A$.

        Remember that $Y$ and $Y_0$ are equipped with topologies
induced by the product topologies on $X$ and $X_0$, respectively.  It
is easy to see that $A$ is continuous as a mapping from ${\bf R}
\times Y_0$ into $Y$, where ${\bf R} \times Y_0$ is equipped with the
product topology associated to the standard topology on ${\bf R}$ and
the topology on $Y_0$ just mentioned.  This uses the fact that $q_l :
{\bf R} \to {\bf R} / r^l \, {\bf Z}$ is continuous for each $l$.  One
can also check that $A$ is a local homeomorphism with respect to these
topologies.  Continuous local inverses for $A$ can be given as in the
previous paragraph, using the fact that $q_0 : {\bf R} \to {\bf R} /
{\bf Z}$ is a local homeomorphism.

\section{A nice mapping, continued}
\label{nice mapping, continued}
\setcounter{equation}{0}

        Consider the metric on ${\bf R} \times Y_0$ defined by
\begin{equation}
\label{D((a, x), (b, y)) = max(|a - b|, rho(x, y))}
        D((a, x), (b, y)) = \max(|a - b|, \rho(x, y)).
\end{equation}
Here $|a|$ is the ordinary absolute value of a real number $a$, so
that $|a - b|$ is the standard metric on the real line, and $\rho(x,
y)$ is the ultrametric on $X_0$ defined in Section \ref{another
metric}.  Thus the topology on ${\bf R} \times Y_0$ determined by
(\ref{D((a, x), (b, y)) = max(|a - b|, rho(x, y))}) is the same as the
product topology associated to the standard topology on ${\bf R}$ and
the usual topology on $Y_0$.  We would like to look more precisely at
the behavior of the mapping $A : {\bf R} \times Y_0 \to Y$ defined in
the previous section with respect to this metric on ${\bf R} \times
Y_0$ and the metric $d(\cdot, \cdot)$ on $Y$ discussed in Section
\ref{nice metric}.

        Note that the derivative of $\exp (i \, t)$ is equal to $i \,
\exp (i \, t)$, which has modulus equal to $1$ for each $t \in {\bf
R}$.  Using this, one can check that
\begin{equation}
\label{|exp (i u) - exp (i v)| le |u - v|}
        |\exp (i \, u) - \exp (i \, v)| \le |u - v|
\end{equation}
for every $u, v \in {\bf R}$, by expressing $\exp (i \, u) - \exp (i
\, v)$ as an integral of $i \, \exp (i \, t)$.  If $|u - v| \le \pi$,
for instance, then we have that
\begin{equation}
\label{|exp (i u) - exp (i v)| ge c_1 |u - v|}
        |\exp (i \, u) - \exp (i \, v)| \ge c_1 \, |u - v|
\end{equation}
for a suitable constant $c_1 > 0$, i.e., $2/\pi$.

        Let $a$ and $b$ be real numbers, and let $\phi_l$, $d_l$, and
$d$ be as in Section \ref{nice metric}.  Thus
\begin{eqnarray}
\label{|phi_l(q_l(a)) - phi_l(q_l(b))| = ... le 2 pi r^{-l} |a - b|}
 |\phi_l(q_l(a)) - \phi_l(q_l(b))| & = & |\exp (2 \, \pi \, i \, r^{-l} \, a)
                                 - \exp (2 \, \pi \, i \, r^{-l} \, b)| \\
                        & \le & 2 \, \pi \, r^{-l} \, |a - b| \nonumber
\end{eqnarray}
for each $l \ge 0$, by (\ref{|exp (i u) - exp (i v)| le |u - v|}).
This implies that
\begin{equation}
\label{d(q(a), q(b)) le 2 pi |a - b|}
        d(q(a), q(b)) \le 2 \, \pi \, |a - b|.
\end{equation}

        Now let $x, y \in Y_0$ be given as well.  If $x = y$, then
\begin{eqnarray}
\label{d(A(a, x), A(b, y)) le 2 pi |a - b| = 2 pi D((a, x), (b, y))}
 d(A(a, x), A(b, y)) & = & d(q(a) + x, q(b) + y) = d(q(a), q(b))  \\
    & \le & 2 \, \pi \, |a - b| = 2 \, \pi \, D((a, x), (b, y)), \nonumber
\end{eqnarray}
using the translation-invariance of $d$ in the second step, and
(\ref{d(q(a), q(b)) le 2 pi |a - b|}) in the third.

        Suppose instead that $x \ne y$, and let $l(x, y)$ be as in
Section \ref{another metric}.  Remember that $x$ and $y$ are
identified with elements of $Y$ by putting $x_0 = y_0 = 0$ in ${\bf R}
/ {\bf Z}$.  If $0 \le j < l(x, y)$, then $x_j = y_j$ by the
definition of $l(x, y)$, and hence
\begin{equation}
\label{d_j(q_j(a) + x_j, q_j(b) + y_j) = d_j(q_j(a), q_j(b))}
        d_j(q_j(a) + x_j, q_j(b) + y_j) = d_j(q_j(a), q_j(b))
\end{equation}
by translation-invariance.  This implies that
\begin{equation}
\label{max_{0 le j < l(x, y)} r^{-j} d_j(q_j(a) + x_j, q_j(b) + y_j) le ...}
 \max_{0 \le j < l(x, y)} r^{-j} \, d_j(q_j(a) + x_j, q_j(b) + y_j)
                                             \le 2 \, \pi \, |a - b|,
\end{equation}
as before.  If $j \ge l(x, y)$, then we have that
\begin{equation}
\label{r^{-j} d_j(q_j(a) + x_j, q_j(b) + y_j) le 2 r^{-l(x, y)}}
        r^{-j} \, d_j(q_j(a) + x_j, q_j(b) + y_j) \le 2 \, r^{-l(x, y)},
\end{equation}
because $d_j \le 2$ automatically, as in (\ref{|phi_l(x_l) -
phi_l(y_l)| le |phi_l(x_l)| + |phi_l(y_l)| = 2}).  This implies that
\begin{equation}
\label{max_{j ge l(x, y)} r^{-j} d_j(q_j(a) + x_j, q_j(b) + y_j) le ...}
        \max_{j \ge l(x, y)} r^{-j} \, d_j(q_j(a) + x_j, q_j(b) + y_j)
                                      \le 2 \, r^{-1} \, \rho(x, y),
\end{equation}
by the definition (\ref{rho(x, y) = r^{-l(x, y) + 1}}) of $\rho(x, y)$.
Combining (\ref{max_{0 le j < l(x, y)} r^{-j} d_j(q_j(a) + x_j, q_j(b)
+ y_j) le ...}) and (\ref{max_{j ge l(x, y)} r^{-j} d_j(q_j(a) + x_j,
q_j(b) + y_j) le ...}), we get that
\begin{equation}
\label{d(A(a, x), A(b, y)) le 2 pi D((a, x), (b, y))}
        d(A(a, x), A(b, y)) \le 2 \, \pi \, D((a, x), (b, y)),
\end{equation}
since $r^{-1} \le 1 \le \pi$.  This also holds when $x = y$, as in
(\ref{d(A(a, x), A(b, y)) le 2 pi |a - b| = 2 pi D((a, x), (b, y))}),
which amounts to taking $l(x, y) = +\infty$ in this argument.

        To get an estimate in the other direction, let us restrict
our attention to $a, b \in {\bf R}$ such that
\begin{equation}
\label{|a - b| le 1/2}
        |a - b| \le 1/2,
\end{equation}
for instance.  Note that
\begin{equation}
\label{d(A(a, x), A(b, y)) ge ... = d_0(q_0(a), q_0(b))}
 d(A(a, x), A(b, y)) \ge d_0(q_0(a) + x_0, q_0(b) + y_0) = d_0(q_0(a), q_0(b)),
\end{equation}
by taking $l = 0$ in the definition (\ref{d(x, y) = max_{l ge 0} r^{-l} 
|phi_l(x_l) - phi_l(y_l)|}) of $d$, and remembering that $x_0 = y_0 = 0$.
Of course,
\begin{eqnarray}
\label{d_0(q_0(a), q_0(b)) = |phi_0(q_0(a)) - phi_0(q_0(b))| = ...}
        d_0(q_0(a), q_0(b)) & = & |\phi_0(q_0(a)) - \phi_0(q_0(b))| \\
  & = & |\exp (2 \, \pi \, i \, a) - \exp (2 \, \pi \, i \, b)|, \nonumber
\end{eqnarray}
so that
\begin{equation}
\label{d(A(a, x), A(b, y)) ge 2 pi c_1 |a - b|}
        d(A(a, x), A(b, y)) \ge 2 \, \pi \, c_1 \, |a - b|
\end{equation}
when $a$, $b$ satisfy (\ref{|a - b| le 1/2}), by (\ref{|exp (i u) -
exp (i v)| ge c_1 |u - v|}).  In particular, we can combine this with
(\ref{d(q(a), q(b)) le 2 pi |a - b|}) to get that
\begin{equation}
\label{d(q(a), q(b)) le c_1^{-1} d(A(a, x), A(b, y))}
        d(q(a), q(b)) \le c_1^{-1} \, d(A(a, x), A(b, y))
\end{equation}
when $a$, $b$ satisfy (\ref{|a - b| le 1/2}).

        Using translation-invariance and then the triangle inequality,
we get that
\begin{equation}
\label{d(x, y) = ... le d(q(a) + x, q(b) + y) + d(q(a), q(b))}
 \quad  d(x, y) = d(q(a) + x, q(a) + y)
                  \le d(q(a) + x, q(b) + y) + d(q(a), q(b)).
\end{equation}
Combining this with (\ref{d(q(a), q(b)) le c_1^{-1} d(A(a, x), A(b, y))}),
it follows that
\begin{equation}
\label{d(x, y) le d(A(a, x), A(b, y)) + c_1^{-1} d(A(a, x), A(b, y))}
        d(x, y) \le d(A(a, x), A(b, y)) + c_1^{-1} \, d(A(a, x), A(b, y))
\end{equation}
when $a$, $b$ satisfy (\ref{|a - b| le 1/2}).  We also know that
$\rho(x, y)$ is bounded by a constant times $d(x, y)$, as in
(\ref{d(x, y) ge r^{-1} |exp (2 pi i r^{-1}) - 1| rho(x, y)}), so that
\begin{equation}
\label{rho(x, y) le constant times d(A(a, x), A(b, y))}
        \rho(x, y) \le r \, |\exp (2 \, \pi \, i \, r^{-1}) - 1|^{-1}
                          \, (1 + c_1^{-1}) \, d(A(a, x), A(b, y))
\end{equation}
when $a$, $b$ satisfy (\ref{|a - b| le 1/2}).  This together with
(\ref{d(A(a, x), A(b, y)) ge 2 pi c_1 |a - b|}) shows that
\begin{equation}
\label{D((a, x), (b, y)) le c_2(r) d(A(a, x), A(b, y))}
        D((a, x), (b, y)) \le c_2(r) \, d(A(a, x), A(b, y))
\end{equation}
when $a$, $b$ satisfy (\ref{|a - b| le 1/2}), where $c_2(r)$ is a
positive real number that depends only on $r$.

        Note that the comparison between $d(A(a, x), A(b, y))$ and
$D((a, x), (b, y))$ in this section would be a bit simpler if we
replaced $\rho(x, y)$ in (\ref{D((a, x), (b, y)) = max(|a - b|, rho(x,
y))}) with $r^{-1} \, \rho(x, y)$, to get the metric
\begin{equation}
\label{D'((a, x), (b, y)) = max(|a - b|, r^{-1} rho(x, y))}
        D'((a, x), (b, y)) = \max(|a - b|, r^{-1} \, \rho(x, y))
\end{equation}
on ${\bf R} \times Y_0$.  Similarly, the comparison between $d(x, y)$
and $\rho(x, y)$ for $x, y \in Y_0$ in Section \ref{another metric}
may be considered as a better comparison between $d(x, y)$ and $r^{-1}
\, \rho(x, y)$.  However, the original definition $\rho(x, y)$ has the
advantage that it corresponds exactly to the $r$-adic metric on ${\bf
Z}_r$, as in Section \ref{r-adic integers}.

\section{Haar measure on $Y$}
\label{haar measure on Y}
\setcounter{equation}{0}

        It is well known that every locally compact commutative
topological group has a nonnegative Borel measure which is invariant
under translations, finite on compact sets, and positive on nonempty
open sets, known as \emph{Haar measure}.  This measure is unique up to
multiplication by a positive real number, at least under some
additional regularity conditions, which are not necessary for the
groups under consideration here.  Of course, Lebesgue measure
satisfies the requirements of Haar measure on the real line as a
locally compact commutative topological group with respect to
addition, and similarly for the unit circle.  Alternatively, one can
start with a nonnegative linear functional on the space of continuous
real or complex-valued functions with compact support on the group
which is invariant under translations and strictly positive for
nonnegative continuous functions that are positive somewhere on the group.
The Riesz representation theorem then leads to a nonnegative Borel measure
on the group with the required properties.

        Let us begin with $Y_0$, which we have identified with a
closed subgroup of $Y$, and which is isomorphic as a topological group
to the group ${\bf Z}_r$ of $r$-adic integers with respect to
addition, as in Section \ref{r-adic integers}.  If we normalize Haar
measure on ${\bf Z}_r$ so that the measure of ${\bf Z}_r$ is equal to
$1$, then it is easy to see that the measure of $r^l \, {\bf Z}_r$ has
to be equal to $r^{-l}$ for each nonnegative integer $l$.  This is
because ${\bf Z}_r / r^l \, {\bf Z}_r$ is isomorphic to ${\bf Z} / r^l
\, {\bf Z}$, so that ${\bf Z}_r$ is the union of $r^l$
pairwise-disjoint translates of $r^l \, {\bf Z}_r$.  One can also
define the Haar integral of a continuous function on ${\bf Z}_r$
directly as a limit of Riemann sums, using this partition of ${\bf
Z}_r$ into translates of $r^l \, {\bf Z}_r$ for each $l \ge 0$.

        If $f$ is a continuous real or complex-valued function on $Y$,
then one can first integrate $f$ over $Y_0$ and its translates in $Y$,
to get a continuous function $f_0$ on $Y$ that is constant on $Y_0$
and its translates in $Y$.  Thus $f_0$ is basically the same as a
continuous function on the unit circle, which can be integrated over
${\bf T}$ in the usual way.  It is easy to see that translations of
$f$ on $Y$ correspond to translations of the function on ${\bf T}$
associated to $f_0$ in a simple way, so that this defines a
translation-invariant integral of continuous functions on $Y$ with the
appropriate positivity properties.  Equivalently, one can average $f$
over translates of the subgroup of $Y_0 \cong {\bf Z}_r$ that
corresponds to $r^l \, {\bf Z}_r$ for any nonnegative integer $l$, to
get a continuous function $f_l$ on $Y$ that is basically the same as a
continuous function on ${\bf R} / r^l \, {\bf Z}$.  One can then take
the average of the resulting function on ${\bf R} / r^l \, {\bf Z}$ to
get a translation-invariant average of $f$ on $Y$ that does not depend
on $l$.

        Note that Haar measure on ${\bf Z}_r$ is Ahlfors regular of
dimension $1$ with respect to the $r$-adic metric on ${\bf Z}_r$, in
the sense that the measure of a ball of radius $t > 0$ is bounded from
above and below by constant multiples of $t$, at least when $t$ is
less than or equal to the diameter of ${\bf Z}_r$, which is $1$.  More
precisely, the closed balls in ${\bf Z}_r$ of radius $r^{-l}$ are the
same as the translates of $r^l \, {\bf Z}_r$, which have measure equal
to $r^{-l}$.  Of course, Lebesgue measure on the real line is also
Ahlfors regular of dimension $1$ with respect to the standard metric
on ${\bf R}$, and Haar measure on the unit circle is Ahlfors regular
of dimension $1$ with respect to the standard metric on ${\bf T}$ as
well.  Similarly, one can check that Haar measure on $Y$ is
Ahlfors-regular of dimension $2$ with respect to the metric $d$
defined in Section \ref{nice metric}.

\section{Continuous functions on $Y$}
\label{continuous functions on Y}
\setcounter{equation}{0}

        Let $\pi_l$ be the $l$th coordinate projection from $X$ onto
${\bf R} / r^l \, {\bf Z}$ for each nonnegative integer $l$, as in
Section \ref{cartesian product}.  Thus the restriction of $\pi_l$ to
$Y$ defines a continuous homomorphism from $Y$ onto ${\bf R} / r^l \,
{\bf Z}$ for each $l$.  In particular, if $g$ is a continuous real or
complex-valued function on ${\bf R} / r^l \, {\bf Z}$, then the
restriction of $g \circ \pi_l$ to $Y$ is a continuous function on $Y$.
These are the same as the continuous functions on $Y$ that are
constant on the translates of the subgroup of $Y_0 \cong {\bf Z}_r$
corresponding to $r^l \, {\bf Z}_r$.

        If $f$ is any continuous real or complex-valued function on
$Y$, then $f$ can be approximated uniformly on $Y$ by functions of
this type, as $l \to \infty$.  One way to see this is to average $f$
over the translates of the subgroups of $Y_0$ corresponds to $r^l \,
{\bf Z}_r$ in ${\bf Z}_r$, with respect to Haar measure on ${\bf
Z}_r$.  These averages will converge uniformly to $f$ as $l \to
\infty$, because of the uniform continuity of $f$ on $Y$, and since
$Y$ is a compact metric space.  Alternatively, $f(q(a))$ is a
continuous function on the real line such that
\begin{equation}
\label{lim_{l to infty} f(q(r^l)) = f(q(0))}
        \lim_{l \to \infty} f(q(r^l)) = f(q(0)),
\end{equation}
since $q(r^l) \to q(0)$ as $l \to \infty$ in $Y$.  This permits one to
approximate the restriction of $f(q(a))$ to $[0, r^l]$ by a continuous
periodic function with period $r^l$.  The latter corresponds exactly
to a continuous function on ${\bf R} / r^l \, {\bf Z}$, whose
composition with $\pi_l$ defines a continuous function on $Y$ as
before.  One can again use the uniform continuity of $f$ to show that
$f$ is uniformly approximated on $Y$ by functions like these as $l \to
\infty$.

        Remember that a \emph{character} on a locally compact
commutative topological group is a continuous homomorphism from that
group into the unit circle ${\bf T}$, as a group with respect to
multiplication of complex numbers.  It is well known that the
characters on ${\bf T}$ are the mappings of the form $z \mapsto z^n$,
where $n$ is an integer.  Equivalently, the characters on ${\bf R} /
r^l \, {\bf Z}$ are given by the integer powers of the mappings
$\phi_l$ defined in Section \ref{nice metric}.  The composition of any
character on ${\bf R} / r^l \, {\bf Z}$ with $\pi_l$ defines a
character on $Y$, since the restriction of $\pi_l$ to $Y$ is a
continuous homomorphism from $Y$ onto ${\bf R} / r^l \, {\bf Z}$.  If
$k$ is an integer greater than or equal to $l$, then $r^k \, {\bf Z}$
is a subgroup of $r^l \, {\bf Z}$, which leads to a continuous
homomorphism from ${\bf R} / r^k \, {\bf Z}$ onto ${\bf R} / r^l \,
{\bf Z}$.  The composition of a character on ${\bf R} / r^l \, {\bf
Z}$ with this homomorphism leads to a character on ${\bf R} / r^k \,
{\bf Z}$, and then to a character on $Y$ by composition with $\pi_k$.
Thus the characters on $Y$ coming from those on ${\bf R} / r^l \, {\bf
Z}$ are contained in the characters on $Y$ coming from those on ${\bf
R} / r^k \, {\bf Z}$ when $l \le k$.

        Conversely, every character on $Y$ comes from one on ${\bf R}
/ r^l \, {\bf Z}$ in this way.  To see this, it suffices to check that
every character on $Y$ is constant on one of the subgroups of $Y_0
\cong {\bf Z}_r$ corresponding to $r^l \, {\bf Z}_r$ for some $l$.
Note that every neighborhood of the additive identity element $0$ in
$Y$ contains these subgroups for sufficiently large $l$.  A character
on $Y$ maps small neighborhoods of $0$ in $Y$ to small neighborhoods
of $1$ in ${\bf T}$, and hence maps these subgroups into small
neighborhoods of $1$ in ${\bf T}$ when $l$ is sufficiently large.
However, the trivial subgroup $\{1\}$ of ${\bf T}$ is the only
subgroup contained in a suitable neighborhood of $1$ in ${\bf T}$,
which implies that characters on $Y$ are constant on these subgroups
when $l$ is sufficiently large, as desired.

\section{Concluding remarks}
\label{concluding remarks}
\setcounter{equation}{0}

        The solenoid $Y$ seems to be an interesting example of a
somewhat exotic ``space of homogeneous type'', in the sense of
\cite{c-w-1, c-w-2}.  Of course, the local geometry of $Y$ is
essentially that of a product of an interval with a Cantor set, but
the global structure is more complicated, since $Y$ is connected in
particular.  In addition, $Y$ has the structure of a compact
commutative topological group, and the geometry on $Y$ is compatible
with this.  It should also be mentioned that for some questions in
analysis, one is probably better off looking at $Y$ as a topological
group, without using this type of geometry.  More precisely, one can
approximate $Y$ by ${\bf R} / r^l \, {\bf Z}$, and functions on $Y$ be
functions on ${\bf R} / r^l \, {\bf Z}$, as in the previous section.

        By way of comparison, one might consider the ordinary product
$W = {\bf T} \times {\bf Z}_r$ of the unit circle and the $r$-adic
integers.  This is also a compact commutative topological group, where
the group operations are defined coordinatewise, and one can get a
natural translation-invariant metric on $W$ by taking the maximum of
the usual metrics on ${\bf T}$ and ${\bf Z}_r$ in their respective
coordinates.  Note that characters on $W$ are given by products of
characters on ${\bf T}$ and ${\bf Z}_r$.  In this case, Haar measure
on $W$ is given by the product of the Haar measures on ${\bf T}$ and
${\bf Z}_r$, and is Ahlfors regular of dimension $2$ in particular.

        Let $M$ be a metric space which is the product of a closed
interval in the real line with the standard metric and another metric
space which is Ahlfors regular of some positive dimension.  As in
Theorem 4.12 in Section 4.4 of \cite{s2}, one can use arguments like
those in \cite{gh} to show that metric doubling measures are
absolutely continuous, with density given by an $A_\infty$ weight.
This is basically the same as absolute continuity properties of
quasisymmetric mappings from $M$ into another Ahlfors regular metric
space of the same dimension, which is the Hausdorff dimension.  This
type of argument is essentially local, and hence works as well for
spaces like $Y$.  Of course, the global structure of $Y$ is important
for the global behavior of quasisymmetric mappings on $Y$ too.

\part{More complicated versions}
\label{more complicated versions}

\section{The $r$-adic metric}
\label{r-adic metric}
\setcounter{equation}{0}

        Let $r = \{r_j\}_{j = 1}^\infty$ be a sequence of integers
with $r_j \ge 2$ for each $j$, and put
\begin{equation}
\label{R_l = prod_{j = 1}^l r_j}
        R_l = \prod_{j = 1}^l r_j
\end{equation}
when $l \ge 1$, and $R_0 = 1$.  In particular, if $r$ is a constant
sequence, so that $r_j = r_1$ for each $j$, then
\begin{equation}
\label{R_l = r_1^l}
        R_l = r_1^l
\end{equation}
for each $l$.  Note that $R_l \ge 2^l$ for each $l \ge 0$, and hence
that $R_l \to \infty$ as $l \to \infty$.  If $r_j \le C$ for some $C
\ge 2$ and each $j \ge 1$, then $R_l \le C^l$ for each $l \ge 0$.

        Let ${\bf Z}$ be the ring of integers, as usual.  If $x \in
{\bf Z}$ and $x \ne 0$, then let $l(x)$ be the largest nonnegative
integer $l$ such that $x$ is an integer multiple of $R_l$, and put
$l(0) = +\infty$.  Observe that
\begin{equation}
\label{l(x + y) ge min(l(x), l(y))}
        l(x + y) \ge \min(l(x), l(y))
\end{equation}
and
\begin{equation}
\label{l(x y) ge max(l(x), l(y))}
        l(x \, y) \ge \max(l(x), l(y))
\end{equation}
for every $x, y \in {\bf Z}$.  If $r$ is a constant sequence, then we
get that
\begin{equation}
\label{l(x y) ge l(x) + l(y)}
        l(x \, y) \ge l(x) + l(y)
\end{equation}
for every $x, y \in {\bf Z}$.  If $r$ is a constant sequence and $r_1$
is a prime number, then
\begin{equation}
\label{l(x y) = l(x) + l(y)}
        l(x \, y) = l(x) + l(y)
\end{equation}
for every $x, y \in {\bf Z}$.

        The \emph{$r$-adic absolute value} of $x \in {\bf Z}$ is
defined by
\begin{equation}
\label{|x|_r = 1/R_{l(x)}}
        |x|_r = 1/R_{l(x)}
\end{equation}
when $x \ne 0$, and $|0|_r = 0$.  Thus
\begin{equation}
\label{|x + y|_r le max(|x|_r, |y|_r)}
        |x + y|_r \le \max(|x|_r, |y|_r)
\end{equation}
and
\begin{equation}
\label{|x y|_r le min(|x|_r, |y|_r)}
        |x \, y|_r \le \min(|x|_r, |y|_r)
\end{equation}
for every $x, y \in {\bf Z}$, by (\ref{l(x + y) ge min(l(x), l(y))})
and (\ref{l(x y) ge max(l(x), l(y))}).  If $r$ is a constant sequence,
then
\begin{equation}
\label{|x y|_r le |x|_r |y|_r}
        |x \, y|_r \le |x|_r \, |y|_r
\end{equation}
for every $x, y \in {\bf Z}$, by (\ref{l(x y) ge l(x) + l(y)}).
Similarly, if $r$ is a constant sequence and $r_1$ is a prime number,
then
\begin{equation}
        |x \, y|_r = |x|_r \, |y|_r
\end{equation}
for every $x, y \in {\bf Z}$, by (\ref{l(x y) = l(x) + l(y)}).  In
this case, $|x|_r$ is the same as the usual $p$-adic absolute value of
$x$, with $p = r_1$.

        The \emph{$r$-adic metric} on ${\bf Z}$ is defined by
\begin{equation}
\label{delta_r(x, y) = |x - y|_r}
        \delta_r(x, y) = |x - y|_r.
\end{equation}
It is easy to see that this satisfies the requirements of a metric on
${\bf Z}$, and is in fact an \emph{ultrametric} on ${\bf Z}$, since
\begin{equation}
\label{delta_r(x, z) le max (delta_r(x, y), delta_r(y, z))}
        \delta_r(x, z) \le \max (\delta_r(x, y), \delta_r(y, z))
\end{equation}
for every $x, y, z \in {\bf Z}$, by (\ref{|x + y|_r le max(|x|_r, |y|_r)}).
By construction, this ultrametric is invariant under translations
on ${\bf Z}$, in the sense that
\begin{equation}
\label{delta_r(x + z, y + z) = delta_r(x, y)}
        \delta_r(x + z, y + z) = \delta_r(x, y)
\end{equation}
for every $x, y, z \in {\bf Z}$.  One can check that addition and
multiplication on ${\bf Z}$ define continuous mappings from ${\bf Z}
\times {\bf Z}$ into ${\bf Z}$ with respect to the topology associated
to this ultrametric, using the corresponding product topology on ${\bf
Z} \times {\bf Z}$.  More precisely, this follows from (\ref{|x + y|_r
le max(|x|_r, |y|_r)}), (\ref{|x y|_r le min(|x|_r, |y|_r)}), and
standard arguments.

\section{$r$-Adic integers}
\label{r-adic integers, 2}
\setcounter{equation}{0}

        As usual, one can take the completion of ${\bf Z}$ as a metric
space with the ultrametric $\delta_r(x, y)$ to get the $r$-adic integers
${\bf Z}_r$.  If $r$ is a constant sequence and $r_1 = p$ is a prime
number, then this reduces to the usual construction of the $p$-adic
integers ${\bf Z}_p$.  For any $r$, addition and multiplication can be
extended to ${\bf Z}_r$, as well as the $r$-adic absolute value and
metric, and with the same properties as before.  In particular, ${\bf
Z}_r$ is a commutative topological ring with respect to the topology
determined by the extension of the $r$-adic metric.  It is easy to see
that ${\bf Z}$ and hence ${\bf Z}_r$ are totally bounded with respect
to the $r$-adic metric, which implies that ${\bf Z}_r$ is compact,
since it is complete.

        Alternatively, consider the Cartesian product
\begin{equation}
\label{X_0 = prod_{l = 1}^infty ({bf Z} / R_l {bf Z})}
        X_0 = \prod_{l = 1}^\infty ({\bf Z} / R_l \, {\bf Z}).
\end{equation}
Thus the elements of $X_0$ are sequences $x = \{x_l\}_{l = 1}^\infty$
with $x_l$ in the quotient ring ${\bf Z} / R_l \, {\bf Z}$ for each
$l$.  Addition and multiplication of elements of $X_0$ can be defined
coordinatewise, so that $X_0$ is a commutative ring.  Using the
product topology on $X_0$ corresponding to the discrete topology on
${\bf Z} / R_l \, {\bf Z}$ for each $l$, $X_0$ becomes a compact
Hausdorff topological space as well.  It is easy to see that addition
and multiplication on $X_0$ define continuous mappings from $X_0
\times X_0$ into $X_0$, using the product topology on $X_0 \times
X_0$, so that $X_0$ is a topological ring.

        If $x, y \in X_0$ and $x \ne y$, then let $l(x, y)$ be the
smallest positive integer $l$ such that $x_l \ne y_l$.  Equivalently,
$l(x, y) - 1$ is the largest nonnegative integer such that $x_l
= y_l$ when $l \le l(x, y) - 1$.  Put
\begin{equation}
\label{rho(x, y) = 1/R_{l(x, y) - 1}}
        \rho(x, y) = 1/R_{l(x, y) - 1}.
\end{equation}
If $x = y$, then we put $l(x, y) = +\infty$ and $\rho(x, y) = 0$.  Thus
\begin{equation}
\label{l(x, y) = l(y, x), 2}
        l(x, y) = l(y, x)
\end{equation}
and
\begin{equation}
\label{rho(x, y) = rho(y, x), 2}
        \rho(x, y) = \rho(y, x)
\end{equation}
for every $x, y \in X_0$.  One can also check that
\begin{equation}
\label{l(x, z) ge min(l(x, y), l(y, z)), 2}
        l(x, z) \ge \min(l(x, y), l(y, z))
\end{equation}
for every $x, y, z \in X_0$, so that
\begin{equation}
\label{rho(x, z) le max(rho(x, y), rho(y, z)), 2}
        \rho(x, z) \le \max(\rho(x, y), \rho(y, z)).
\end{equation}
This implies that $\rho(x, y)$ defines an ultrametric on $X_0$, and it
is easy to see that the topology on $X_0$ corresponding to $\rho(x,
y)$ is the same as the product topology mentioned in the previous
paragraph.  Note that
\begin{equation}
\label{l(x + z, y + z) = l(x, y), 2}
        l(x + z, y + z) = l(x, y)
\end{equation}
for every $x, y, z \in X_0$, and hence that
\begin{equation}
\label{rho(x + z, y + z) = rho(x, y), 2}
        \rho(x + z, y + z) = \rho(x, y),
\end{equation}
so that $\rho(x, y)$ is invariant under translations on $X_0$.
Because $X_0$ is compact, it is complete as a metric space, which can
also be verified directly from the definitions.

        Let $\widetilde{q}_l$ be the standard quotient map from ${\bf
Z}$ onto ${\bf Z} / R_l \, {\bf Z}$ for each $l \ge 1$, which is a
ring homomorphism.  Put
\begin{equation}
\label{widetilde{q}(a) = {widetilde{q}_l(a)}_{l = 1}^infty, 2}
        \widetilde{q}(a) = \{\widetilde{q}_l(a)\}_{l = 1}^\infty
\end{equation}
for each $a \in {\bf Z}$, which is a ring homomorphism from ${\bf Z}$
into $X_0$ with trivial kernel.  If $a, b \in {\bf Z}$, then it is
easy to see that
\begin{equation}
\label{l(widetilde{q}(a), widetilde{q}(b)) - 1 = l(a - b), 2}
        l(\widetilde{q}(a), \widetilde{q}(b)) - 1 = l(a - b),
\end{equation}
and hence that
\begin{equation}
\label{rho(widetilde{q}(a), widetilde{q}(b)) = |a - b|_r = delta_r(a, b)}
        \rho(\widetilde{q}(a), \widetilde{q}(b)) = |a - b|_r = \delta_r(a, b).
\end{equation}
Thus $\widetilde{q}$ is an isometric embedding of ${\bf Z}$ with the
$r$-adic metric into $X_0$ with the ultrametric $\rho(x, y)$, which
implies that the completion ${\bf Z}_r$ of ${\bf Z}$ with respect to
the $r$-adic metric can be identified with the closure of
$\widetilde{q}({\bf Z})$ in $X_0$, since $X_0$ is complete.

        As usual, there is a natural ring homomorphism from ${\bf Z} /
R_{l + 1} \, {\bf Z}$ onto ${\bf Z} / R_l \, {\bf Z}$ for each $l \ge
0$, because $R_{l + 1} \, {\bf Z} \subseteq R_l \, {\bf Z}$.  An
element $x = \{x_l\}_{l = 1}^\infty$ of $X_0$ is said to be a
\emph{coherent sequence} if $x_l$ is the image in ${\bf Z} / R_l \,
{\bf Z}$ of $x_{l + 1} \in {\bf Z} / R_{l + 1} \, {\bf Z}$ for each $l
\ge 0$.  In particular, $\widetilde{q}(a)$ is a coherent sequence for
each $a \in {\bf Z}$, because $\widetilde{q}_l$ is the same as the
composition of $\widetilde{q_{l + 1}}$ with the natural homomorphism
from ${\bf Z} / R_{l + 1} \, {\bf Z}$ onto ${\bf Z} / R_l \, {\bf Z}$
for each $l$.  The set $Y_0$ of coherent sequences is a closed subring
of $X_0$, and one can check that $\widetilde{q}({\bf Z})$ is dense in
$Y_0$, so that $Y_0$ is the same as the closure in $X_0$ of
$\widetilde{q}({\bf Z})$.  Thus the completion ${\bf Z}_r$ of ${\bf
Z}$ with respect to the $r$-adic metric can be identified with $Y_0$.

        Let $n$ be a positive integer, and let $Y_n$ be the set of
coherent sequences $x = \{x_l\}_{l = 1}^\infty$ such that $x_n = 0$ in
${\bf Z} / R_n \, {\bf Z}$.  This implies that $x_l = 0$ in ${\bf Z} /
R_l \, {\bf Z}$ when $l \le n$, because of coherence.  It is easy to
see that $Y_n$ is a closed subring of $X_0$ which is an ideal in
$Y_0$.  If $a \in {\bf Z}$, then $\widetilde{q}(a) \in Y_n$ if and
only if $a \in R_n \, {\bf Z}$, and $Y_n$ is the same as the closure
in $X_0$ of $\widetilde{q}(R_n \, {\bf Z})$.

        Equivalently, $Y_n$ is the kernel of the homomorphism from
$Y_0$ into ${\bf Z} / R_n \, {\bf Z}$ that sends a coherent sequence
$x = \{x_l\}_{l = 1}^\infty$ to its $n$th term $x_n$.  More precisely,
this is a homomorphism from $Y_0$ onto ${\bf Z} / R_n \, {\bf Z}$,
because its composition with $\widetilde{q} : {\bf Z} \to Y_0$ is the
quotient homomorphism $\widetilde{q}_n$ from ${\bf Z}$ onto ${\bf Z} /
R_n \, {\bf Z}$.  Thus the quotient $Y_0 / Y_n$ is isomorphic as a
commutative ring to ${\bf Z} / R_n \, {\bf Z}$.

\section{Haar measure on $Y_0 \cong {\bf Z}_r$}
\label{haar measure on Y_0 cong Z_r}
\setcounter{equation}{0}

        Let $\mu_0$ be Haar measure on $Y_0$, normalized so that
\begin{equation}
\label{mu_0(Y_0) = 1}
        \mu_0(Y_0) = 1.
\end{equation}
Note that $Y_n$ is both relatively open and closed in $Y_0$ for each
$n \ge 1$.  It is easy to see that
\begin{equation}
\label{mu_0(Y_n) = 1/R_n}
        \mu_0(Y_n) = 1/R_n
\end{equation}
for each $n \ge 1$, because $Y_0 / Y_n \cong {\bf Z} / R_n \, {\bf Z}$,
so that $Y_0$ is the union of $R_n$ pairwise-disjoint translates of $Y_n$.

        With respect to the restriction of the ultrametric $\rho(x,
y)$ on $X_0$ to $Y_0$, $Y_n$ is the same as the closed ball in $Y_0$
centered at $0$ and with radius $1/R_n$, and every closed ball in
$Y_0$ with radius $1/R_n$ is a translate of $Y_n$.  Thus
(\ref{mu_0(Y_n) = 1/R_n}) may be considered as a very precise form of
Ahlfors regularity of dimension $1$ for radii of the form $1/R_n$.  In
particular, if the original sequence of $r_j$'s bounded, then it is
easy to see that $Y_0$ is Ahlfors regular of dimension $1$.  Of
course, the $r_j$'s are bounded when they are all equal to each other,
in which case $Y_0$ enjoys additional self-similarity properties.
However, if the $r_j$'s are not bounded, then $\mu_0$ is not even a
doubling measure on $Y_0$, and $Y_0$ does not satisfy a doubling
condition as a metric space.

        Even if the $r_j$'s are not bounded, the fact that the metric
on $Y_0$ is an ultrametric implies that any two balls in $Y_0$ are
either disjoint, or one of the balls is contained in the other.  Given
any collection of balls in $Y_0$, one can take the maximal balls in
the collection to get a sub-collection of pairwise-disjoint balls with
the same union.  One can also look at this in terms of martingales,
using the partitions of $Y_0$ obtained from the translations of $Y_n$
for each $n$.  Thus one can get the usual estimates for the
Hardy--Littlewood maximal function on $Y_0$, for instance, even when
the $r_j$'s are not bounded.

        If $H^1(E)$ denotes the one-dimensional Hausdorff measure of a
set $E \subseteq Y_0$ with respect to the restriction of the
ultrametric $\rho(x, y)$ on $X_0$ to $Y_0$, then it is easy to see that
\begin{equation}
\label{H^1(Y_n) le  1/R_n}
        H^1(Y_n) \le  1/R_n,
\end{equation}
for each $n \ge 0$, by considering coverings of $Y_n$ by translates of
$Y_k$ when $k \ge n$.  One can get the opposite inequality by
comparing other coverings of $Y_n$ with these, so that
\begin{equation}
\label{H^1(Y_n) = 1/R_n}
        H^1(Y_n) = 1/R_n
\end{equation}
for each $n$.  Of course, Hausdorff measure of any dimension is
automatically invariant under translations on $Y_0$, because the
ultrametric $\rho(x, y)$ is invariant under translations.  It follows
that the normalized Haar measure $\mu_0$ on $Y_0$ is the same as
one-dimensional Hausdorff measure on $Y_0$.

\section{Quotients of ${\bf R}$}
\label{quotients of R}
\setcounter{equation}{0}

        Consider the Cartesian product
\begin{equation}
\label{X = prod_{l = 0}^infty ({bf R} / R_l {bf Z})}
        X = \prod_{l = 0}^\infty ({\bf R} / R_l \, {\bf Z}).
\end{equation}
Here ${\bf R} / R_l \, {\bf Z}$ refers to the quotient of ${\bf R}$ as
a commutative group by the subgroup $R_l \, {\bf Z}$, which is a
commutative group as well.  Thus $X$ is a commutative group too, where
the group operations are defined coordinatewise.  Using the usual
quotient topology on ${\bf R} / R_l \, {\bf Z}$ for each $l$, $X$
becomes a compact Hausdorff topological space with respect to the
product topology.  It is easy to see that the group operations on $X$
are continuous with respect to this topology, so that $X$ is a
topological group.  We can identify the Cartesian product $X_0$ in
(\ref{X_0 = prod_{l = 1}^infty ({bf Z} / R_l {bf Z})}) with a subset
of $X$, by extending each $x = \{x_l\}_{l = 1}^\infty$ in $X_0$ to $l
= 0$ by taking $x_0 = 0$ in ${\bf R} / {\bf Z}$.  This is very
natural, since $R_0 = 1$ and hence ${\bf Z} / R_0 \, {\bf Z}$ is the
trivial group.  More precisely, $X_0$ corresponds to a closed subgroup
of $X$ in this way, and the topology induced on $X_0$ by the one on
$X$ is the same as the topology on $X_0$ defined in Section
\ref{r-adic integers, 2}.  The group structure on $X_0$ as a subgroup
of $X$ is the same as the additve group structure on $X_0$ as a
commutative ring, as before.

        Because $R_{l + 1} \, {\bf Z}$ is a subgroup of $R_l \, {\bf
Z}$, there is a natural homomorphism from ${\bf R} / R_{l + 1} \, {\bf
Z}$ onto ${\bf R} / R_l \, {\bf Z}$ for each $l \ge 0$.  An element $x
= \{x_l\}_{l = 0}^\infty$ of $X_0$ is said to be a \emph{coherent
sequence} if $x_l$ is the image of $x_{l + 1} \in {\bf R} / R_{l + 1}
\, {\bf Z}$ in ${\bf R} / R_l \, {\bf Z}$ for each $l$.  As usual, the
set $Y$ of coherent sequences in $X$ is a closed subgroup of $X$.
This coherence condition reduces to the previous one for elements of
$X_0$, so that the set $Y_0$ of coherent sequences in $X_0$
corresponds exactly to the intersection of $X_0$ with $Y$ in $X$.
Note that $Y_0$ corresponds to a closed subgroup of $Y$.

        Let $q_l$ be the usual quotient mapping from ${\bf R}$ onto
${\bf R} / R_l \, {\bf Z}$ for each $l \ge 0$, and put
\begin{equation}
\label{q(a) = {q_l(a)}_{l = 0}^infty, 2}
        q(a) = \{q_l(a)\}_{l = 0}^\infty
\end{equation}
for each $a \in {\bf R}$.  This defines a continuous homomorphism from
${\bf R}$ into $X$ with trivial kernel, whose restriction to ${\bf Z}$
corresponds exactly to the mapping $\widetilde{q}$ in
(\ref{widetilde{q}(a) = {widetilde{q}_l(a)}_{l = 1}^infty, 2}).  As
before, $q(a)$ is a coherent sequence in $X$ for each $a \in {\bf R}$,
because $q_l$ is the same as the composition of $q_{l + 1}$ with the
natural homomorphism from ${\bf R} / R_{l + 1} \, {\bf Z}$ onto ${\bf
R} / R_l \, {\bf Z}$ for each $l$.  Thus $q({\bf R}) \subseteq Y$, and
one can check that $q({\bf R})$ is dense in $Y$, so that the closure
of $q({\bf R})$ in $X$ is equal to $Y$.  This implies that $Y$ is
connected, while $Y_0$ is totally disconnected.

        Let $\pi_n$ be the $n$th coordinate projection of $X$ onto
${\bf R} / R_n \, {\bf Z}$, so that
\begin{equation}
\label{pi_n(x) = x_n}
        \pi_n(x) = x_n
\end{equation}
for each $x = \{x_l\}_{l = 0}^\infty \in X$ and nonnegative integer $n$.
This is a continuous group homomorphism from $X$ onto ${\bf R} / R_n
\, {\bf Z}$ for each $n \ge 0$, and we are especially interested in the
restriction of $\pi_n$ to $Y$.  Note that $\pi_n$ maps $Y$ onto ${\bf R} /
R_n \, {\bf Z}$ for each $n$, because $q({\bf R}) \subseteq Y$ and
$\pi_n \circ q = q_n$ maps ${\bf R}$ onto ${\bf R} / R_n \, {\bf Z}$.
The kernel of the restriction of $\pi_n$ to $Y$ corresponds exactly to
the subgroup $Y_n$ of $Y_0$ defined in Section \ref{r-adic integers, 2}
for each $n \ge 0$.  By construction, the restriction of $\pi_n$ to $Y$
is equal to the composition of the restriction of $\pi_{n + 1}$ to $Y$
with the natural homomorphism from ${\bf R} / R_{n + 1} \, {\bf Z}$
onto ${\bf R} / R_n \, {\bf Z}$ for each $n$.

\section{Metrics on $X$}
\label{metrics on X}
\setcounter{equation}{0}

        Let $\phi_l$ be the standard isomorphism from ${\bf R} / R_l
\, {\bf Z}$ onto the unit circle, so that
\begin{equation}
\label{phi_l(q_l(a)) = exp (2 pi i R_l^{-1} a)}
        \phi_l(q_l(a)) = \exp (2 \, \pi \, i \, R_l^{-1} \, a)
\end{equation}
for each $a \in {\bf R}$.  Thus
\begin{equation}
\label{d_l(x_l, y_l) = |phi_l(x_l) - phi_l(y_l)|, 2}
        d_l(x_l, y_l) = |\phi_l(x_l) - \phi_l(y_l)|
\end{equation}
defines a metric on ${\bf R} / R_l \, {\bf Z}$, which determines the same
topology on ${\bf R} / R_l \, {\bf Z}$ as the usual quotient topology.
Note that this metric is invariant under translations on ${\bf R} / 
R_l \, {\bf Z}$, and that
\begin{equation}
\label{d_l(x_l, y_l) le |phi_l(x_l)| + |phi_l(y_l)| = 2}
        d_l(x_l, y_l) \le |\phi_l(x_l)| + |\phi_l(y_l)| = 2
\end{equation}
for every $x_l, y_l \in {\bf R} / R_l \, {\bf Z}$.  If $x_l$, $y_l$
are distinct elements of ${\bf Z} / R_l \, {\bf Z}$ for some $l \ge
1$, then
\begin{equation}
\label{d_l(x_l, y_l) ge |exp (2 pi i R_l^{-1}) - 1|}
        d_l(x_l, y_l) \ge |\exp (2 \, \pi \, i \, R_l^{-1}) - 1|.
\end{equation}
If in addition the images of $x_l$ and $y_l$ in ${\bf Z} / R_{l - 1}
\, {\bf Z}$ are equal, then we get that
\begin{equation}
\label{d_l(x_l, y_l) ge |exp (2 pi i r_l^{-1}) - 1|}
        d_l(x_l, y_l) \ge |\exp (2 \, \pi \, i \, r_l^{-1}) - 1|,
\end{equation}
which is stronger than (\ref{d_l(x_l, y_l) ge |exp (2 pi i R_l^{-1}) - 1|}).

        Let $t = \{t_l\}_{l = 0}^\infty$ be a sequence of positive
real numbers that converges to $0$ in ${\bf R}$, and put
\begin{equation}
\label{d(x, y) = max_{l ge 0} t_l d_l(x_l, y_l)}
        d(x, y) = \max_{l \ge 0} t_l \, d_l(x_l, y_l)
\end{equation}
for each $x, y \in X$.  It is easy to see that the maximum is always
attained under these conditions, and that $d(x, y)$ defines a
translation-invariant metric on $X$ for which the corresponding
topology is the product topology mentioned earlier.  In particular,
the restriction of $d(x, y)$ to $x, y \in Y$ defines a
translation-invariant metric on $Y$, for which the corresponding
topology is the same as the one induced by the product topology on
$X$.  Similarly, if we identify $X_0$, $Y_0$ with subsets of $X$, then
the restriction of $d(x, y)$ to these subsets determine metrics on
$X_0$, $Y_0$ for which the corresponding topologies are the same as
before.  Let us compare this with the ultrametric $\rho(x, y)$ on
$X_0$ defined in Section \ref{r-adic integers, 2}.

        Let $x = \{x_l\}_{l = 1}^\infty, y = \{y_l\}_{l = 1}^\infty
\in X_0$ be given, which can be identified with elements of $X$ by
taking $x_0 = y_0 = 0$ in ${\bf R} / {\bf Z}$, as usual.  We may as
well suppose that $x \ne y$, since otherwise $d(x, y) = \rho(x, y) =
0$.  If $l(x, y)$ is the smallest positive integer $l$ such that $x_l
\ne y_l$, as in Section \ref{r-adic integers, 2}, then we get that
\begin{equation}
\label{d(x, y) = ... le 2 max_{l ge l(x, y)} t_l}
        d(x, y) = \max_{l \ge l(x, y)} t_l \, d_l(x_l, y_l)
                   \le 2 \, \max_{l \ge l(x, y)} t_l.
\end{equation}
If the $t_l$'s are monotone decreasing, then this reduces to
\begin{equation}
\label{d(x, y) le 2 t_{l(x, y)}}
        d(x, y) \le 2 \, t_{l(x, y)}.
\end{equation}
If we take $t_0 = 1$ and $t_l = 1/R_{l - 1}$ when $l \ge 1$, then we get that
\begin{equation}
\label{d(x, y) le 2 / R_{l(x, y) - 1} = 2 rho(x, y)}
        d(x, y) \le 2 / R_{l(x, y) - 1} = 2 \, \rho(x, y).
\end{equation}
In the case where $r = \{r_j\}_{j = 1}^\infty$ is a constant sequence,
so that $R_l = r_1^l$ for each $l \ge 0$, this is the same as taking
$t_0 = 1$ and $t_l = r_1^{-l + 1}$ when $l \ge 1$.  Although $t_l =
r_1^{-l}$ may be appealing in some ways, this slightly different
choice for $t_l$ has other advantages.

        In the other direction, we can take $l = l(x, y)$ in
(\ref{d(x, y) = max_{l ge 0} t_l d_l(x_l, y_l)}), to get that
\begin{equation}
\label{d(x, y) ge ... ge t_{l(x, y)} |exp (2 pi i R_{l(x, y)}^{-1}) - 1|}
 d(x, y) \ge t_{l(x, y)} \, d_l(x_l, y_l)
          \ge t_{l(x, y)} \, |\exp (2 \, \pi \, i \, R_{l(x, y)}^{-1}) - 1|,
\end{equation}
using (\ref{d_l(x_l, y_l) ge |exp (2 pi i R_l^{-1}) - 1|}) in the
second step.  If $x, y \in Y_0$, then we can use (\ref{d_l(x_l, y_l)
ge |exp (2 pi i r_l^{-1}) - 1|}) instead of (\ref{d_l(x_l, y_l) ge
|exp (2 pi i R_l^{-1}) - 1|}) to get that
\begin{equation}
\label{d(x, y) ge t_{l(x, y)} |exp (2 pi i r_{l(x, y)}^{-1}) - 1|}
 d(x, y) \ge t_{l(x, y)} \, |\exp (2 \, \pi \, i \, r_{l(x, y)}^{-1}) - 1|.
\end{equation}
More precisely, $x_{l(x, y) - 1} = y_{l(x, y) - 1}$ by definition of
$l(x, y)$, which implies that the images of $x_{l(x, y)}$ and $y_{l(x,
y)}$ in ${\bf Z} / R_{l(x, y) - 1} \, {\bf Z}$ are the same when $x, y
\in Y_0$, by coherence.  If we take $t_0 = 1$ and $t_l = 1/R_{l - 1}$
when $l \ge 0$, as before, then (\ref{d(x, y) ge t_{l(x, y)} |exp (2
pi i r_{l(x, y)}^{-1}) - 1|}) becomes
\begin{equation}
\label{d(x, y) ge |exp (2 pi i r_{l(x, y)}^{-1}) - 1| rho(x, y)}
 d(x, y) \ge |\exp (2 \, \pi \, i \, r_{l(x, y)}^{-1}) - 1| \, \rho(x, y).
\end{equation}
This implies that
\begin{equation}
\label{d(x, y) ge c_0 rho(x, y)}
        d(x, y) \ge c_0 \, \rho(x, y)
\end{equation}
for some $c_0 > 0$ and every $x, y \in Y_0$ when the $r_j$'s are
bounded.

\section{A nice mapping}
\label{nice mapping, 2}
\setcounter{equation}{0}

        Consider the mapping $A : {\bf R} \times Y_0 \to Y$ defined by
\begin{equation}
\label{A(a, x) = q(a) + x, 2}
        A(a, x) = q(a) + x,
\end{equation}
where $q : {\bf R} \to Y$ is as in (\ref{q(a) = {q_l(a)}_{l =
0}^infty, 2}), and $x = \{x_l\}_{l = 1}^\infty \in Y_0$ is identified
with an element of $Y$ by setting $x_0 = 0$ in ${\bf R} / {\bf Z}$, as
usual.  More precisely, the sum $q(a) + x$ uses the group structure on
$Y$ as a subgroup of $X$ as in Section \ref{quotients of R}.  Note
that $A$ is a homomorphism from ${\bf R} \times Y_0$ into $Y$ with
respect to coordinatewise addition on ${\bf R} \times Y_0$, because
$q$ is a homomorphism from ${\bf R}$ into $Y$.

        If $(a, x)$ is in the kernel of $A$, then $q(a) + x = 0$ in
$Y$, and hence $q_0(a) = 0$ in ${\bf R} / {\bf Z}$, since $x_0 = 0$ by
construction.  Thus $a \in {\bf Z}$, which implies that $x =
\widetilde{q}(-a)$ in $Y_0$, where $\widetilde{q} : {\bf Z} \to Y_0$
is as in (\ref{widetilde{q}(a) = {widetilde{q}_l(a)}_{l = 1}^infty,
2}).  Conversely, if $a \in {\bf Z}$ and $x = \widetilde{q}(-a)$ in
$Y_0$, then $A(a, x) = 0$.

        Let us check that $A$ maps ${\bf R} \times Y_0$ onto $Y$.  If
$y$ is any element of $Y$, then we can first choose $a \in {\bf R}$
such that $y_0 = q_0(a)$ in ${\bf R} / {\bf Z}$.  Hence the $l = 0$
component of $x = y - q(a)$ is equal to $0$, so that $x$ corresponds
to an element of $Y_0$, and $y = A(a, x)$, as desired.

        It is easy to see that $A$ is continuous as a mapping from
${\bf R} \times Y_0$ into $Y$, using the standard topology on ${\bf
R}$, the topologies already discussed on $Y_0$ and $Y$, and the
corresponding product topology on ${\bf R} \times Y_0$.  One can also
check that $A$ is a local homeomorphism, where local inverses for $A$
can be given in terms of local inverses for $q_0 : {\bf R} \to {\bf R}
/ {\bf Z}$, as in the previous paragraph.

\section{A nice mapping, continued}
\label{nice mapping, continued, 2}
\setcounter{equation}{0}

        Consider the metric on ${\bf R} \times Y_0$ defined by
\begin{equation}
\label{D((a, x), (b, y)) = max (|a - b|, rho(x, y))}
        D((a, x), (b, y)) = \max (|a - b|, \rho(x, y)),
\end{equation}
where $\rho(x, y)$ is the ultrametric defined on $X_0$ as in Section
\ref{r-adic integers, 2}.  Note that this is a translation-invariant
metric on ${\bf R} \times Y_0$, and that the topology on ${\bf R}
\times Y_0$ determined by this metric is the same as the product
topology associated to the standard topology on ${\bf R}$ and the
usual topology on $Y_0$.  Throughout this section, we let $d(x, y)$ be
the metric on $X$ in Section \ref{metrics on X}, with $t_0 = 1$ and
$t_l = 1/R_{l - 1}$ when $l \ge 1$.  We would like to look at the
behavior of the mapping $A$ from the previous section with respect to
(\ref{D((a, x), (b, y)) = max (|a - b|, rho(x, y))}) on ${\bf R}
\times Y_0$ and $d(x, y)$ on $Y$.

        Remember that
\begin{equation}
\label{|exp (i u) - exp (i v)| le |u - v|, 2}
        |\exp (i \, u) - \exp (i \, v)| \le |u - v|
\end{equation}
for every $u, v \in {\bf R}$, and that
\begin{equation}
\label{|exp (i u) - exp (i v)| ge c_1 |u - v|, 2}
        |\exp (i \, u) - \exp (i \, v)| \ge c_1 \, |u - v|
\end{equation}
for a suitable constant $c_1 > 0$ when $|u - v| \le \pi$.  Thus
\begin{eqnarray}
\label{|phi_l(q_l(a)) - phi_l(q_l(b))| = ... le 2 pi R_l^{-1} |a - b|}
        |\phi_l(q_l(a)) - \phi_l(q_l(b))|
             & = & |\exp (2 \, \pi \, i \, R_l^{-1} \, a)
                         - \exp (2 \, \pi \, i \, R_l^{-1} \,  b)| \\
            & \le & 2 \, \pi \, R_l^{-1} \, |a - b|           \nonumber
\end{eqnarray}
for every $a, b \in {\bf R}$, and hence
\begin{equation}
\label{d(q(a), q(b)) le 2 pi |a - b|, 2}
        d(q(a), q(b)) \le 2 \, \pi \, |a - b|.
\end{equation}
This only uses the fact that $t_l \le 1$ for each $l \ge 0$.

        Now let $x, y \in Y_0$ be given.  If $x = y$, then
\begin{equation}
\label{d(A(a, x), A(b, y)) = d(q(a) + x, q(b) + y) = d(q(a), q(b))}
        d(A(a, x), A(b, y)) = d(q(a) + x, q(b) + y) = d(q(a), q(b)),
\end{equation}
by translation-invariance, so that
\begin{equation}
\label{d(A(a, x), A(b, y)) le 2 pi |a - b| = 2 pi D((a, x), (b, y)), 2}
 d(A(a, x), A(b, y)) \le 2 \, \pi \, |a - b| = 2 \, \pi \, D((a, x), (b, y)),
\end{equation}
by (\ref{d(q(a), q(b)) le 2 pi |a - b|, 2}).  Otherwise, suppose that
$x \ne y$, and let $l(x, y)$ be as in Section \ref{r-adic integers, 2}.
Remember that $x = \{x_l\}_{l = 1}^\infty$ and $y = \{y_l\}_{l = 1}^\infty$
are extended to $l = 0$ by putting $x_0 = y_0 = 0$ in ${\bf R} / {\bf Z}$,
and that $x_j = y_j$ when $j < l(x, y)$.  This implies that
\begin{equation}
\label{d_j(q_j(a) + x_j, q_j(b) + y_j) = d_j(q_j(a), q_j(b)), 2}
        d_j(q_j(a) + x_j, q_j(b) + y_j) = d_j(q_j(a), q_j(b))
\end{equation}
when $j < l(x, y)$, by translation-invariance, so that
\begin{equation}
\label{max_{0 le j < l(x, y)} t_j d_j(q_j(a) + x_j, q_j(b) + y_j) le ...}
        \max_{0 \le j < l(x, y)} t_j \, d_j(q_j(a) + x_j, q_j(b) + y_j)
                                                 \le 2 \, \pi \, |a - b|,
\end{equation}
as before.  If $j \ge l(x, y) \ge 1$, then $t_j = 1/R_{j - 1}$, and
\begin{equation}
\label{t_j d_j(q_j(a) + x_j, q_j(b) + y_j) le ... le 2 / R_{l(x, y) - 1}}
        t_j \, d_j(q_j(a) + x_j, q_j(b) + y_j) \le 2 / R_{j - 1}
                                              \le 2 / R_{l(x, y) - 1},
\end{equation}
since $d_j \le 2$ automatically.  Thus
\begin{equation}
\label{max_{j ge l(x, y)} t_j d_j(q_j(a) + x_j, q_j(b) + y_j) le 2rho(x, y), 2}
        \max_{j \ge l(x, y)} t_j \, d_j(q_j(a) + x_j, q_j(b) + y_j)
                                                   \le 2 \, \rho(x, y),
\end{equation}
and hence
\begin{equation}
\label{d(A(a, x), A(b, y)) le 2 pi D((a, x), (b, y)), 2}
        d(A(a, x), A(b, y)) \le 2 \, \pi \, D((a, x), (b, y)).
\end{equation}

        To get an estimate in the other direction, we restrict our
attention to $a, b \in {\bf R}$ that satisfy
\begin{equation}
\label{|a - b| le 1/2, 2}
        |a - b| \le 1/2,
\end{equation}
for instance.  Observe that
\begin{equation}
\label{d(A(a, x), A(b, y)) ge ... = d_0(q_0(a), q_0(b)), 2}
 \quad  d(A(a, x), A(b, y)) \ge d_0(q_0(a) + x_0, q_0(b) + y_0)
                                   = d_0(q_0(a), q_0(b)),
\end{equation}
by taking $l = 0$ in the definition of $d$, and remembering that $t_0 = 1$
and $x_0 = y_0 = 0$.  By the definitions of $d_0$ and $\phi_0$, we have that
\begin{equation}
\label{d_0(q_0(a), q_0(b)) = |exp (2 pi i a) - exp (2 pi i b)|}
 d_0(q_0(a), q_0(b)) = |\exp (2 \, \pi \, i \, a) - \exp (2 \, \pi \, i \, b)|,
\end{equation}
and hence that
\begin{equation}
\label{d(A(a, x), A(b, y)) ge 2 pi c_1 |a - b|, 2}
        d(A(a, x), A(b, y)) \ge 2 \, \pi \, c_1 \, |a - b|,
\end{equation}
when $a, b \in {\bf R}$ satisfy (\ref{|a - b| le 1/2, 2}), because of
(\ref{|exp (i u) - exp (i v)| ge c_1 |u - v|, 2}).  Combining this
with (\ref{d(q(a), q(b)) le 2 pi |a - b|, 2}), we get that
\begin{equation}
\label{d(q(a), q(b)) le c_1^{-1} d(A(a, x), A(b, y)), 2}
        d(q(a), q(b)) \le c_1^{-1} \, d(A(a, x), A(b, y))
\end{equation}
when $a, b \in {\bf R}$ satisfy (\ref{|a - b| le 1/2, 2}).  We also have that
\begin{equation}
\label{d(x, y) = ... le d(q(a) + x, q(b) + y) + d(q(a), q(b)), 2}
  \qquad d(x, y) = d(q(a) + x, q(a) + y)
               \le d(q(a) + x, q(b) + y) + d(q(a), q(b)),
\end{equation}
by translation-invariance and the triangle inequality, so that
\begin{equation}
\label{d(x, y) le (1 + c_1^{-1}) d(A(a, x), A(b, y))}
        d(x, y) \le (1 + c_1^{-1}) \, d(A(a, x), A(b, y))
\end{equation}
when $a, b \in {\bf R}$ satisfy (\ref{|a - b| le 1/2, 2}).  
If the $r_j$'s are bounded, then (\ref{d(x, y) ge c_0 rho(x, y)})
holds for some $c_0 > 0$, and we get that
\begin{equation}
\label{c_0 rho(x, y) le (1 + c_1^{-1}) d(A(a, x), A(b, y))}
        c_0 \, \rho(x, y) \le (1 + c_1^{-1}) \, d(A(a, x), A(b, y))
\end{equation}
when $a, b \in {\bf R}$ satisfy (\ref{|a - b| le 1/2, 2}).  Combining
this with (\ref{d(A(a, x), A(b, y)) ge 2 pi c_1 |a - b|, 2}), we get
that
\begin{equation}
\label{D((a, x), (b, y)) le c_2 d(A(a, x), A(b, y))}
        D((a, x), (b, y)) \le c_2 \, d(A(a, x), A(b, y))
\end{equation}
for a suitable constant $c_2$ when $a, b \in {\bf R}$ satisfy
(\ref{|a - b| le 1/2, 2}) and the $r_j$'s are bounded.

\section{Another metric on $Y$}
\label{another metric on Y}
\setcounter{equation}{0}

        Even if the $r_j$'s are not bounded, we can simply choose a
metric on $Y$ which approximates $D$ on ${\bf R} \times Y_0$.  The
easiest way to do that is to take
\begin{eqnarray}
\label{Delta(u, v) = inf {D(A(a, x), A(b, y)) : A(a, x) = u, A(b, y) = v}}
        \Delta(u, v) & = & \inf \{D(A(a, x), A(b, y)) : 
                              (a, x), (b, y) \in {\bf R} \times Y_0, \\
  & & \qquad\qquad\qquad\qquad\qquad \,\ A(a, x) = u, A(b, y) = v\} \nonumber
\end{eqnarray}
for each $u, v \in Y$.  Note that
\begin{equation}
\label{Delta(u, v) = Delta(v, u)}
        \Delta(u, v) = \Delta(v, u)
\end{equation}
for every $u, v \in Y$, and that $\Delta(u, v)$ is
translation-invariant on $Y$, because $D$ is symmetric and
translation-invariant on ${\bf R} \times Y_0$.  If $(a, x) \in {\bf R}
\times Y_0$ satisfies $A(a, x) = u$, then
\begin{equation}
\label{Delta(u, v) = inf {D(A(a, x), A(b, y)) : A(b, y) = v}}
 \Delta(u, v) = \inf \{D(A(a, x), A(b, y)) : (b, y) \in {\bf R} \times Y_0,
                                                \ A(b, y) = v\}
\end{equation}
for every $v \in Y$, since one can use the translation-invariance of
$D$ on ${\bf R} \times Y_0$ to simultaneously translate
representatives of $u$ and $v$ in ${\bf R} \times Y_0$ to reduce to
the case where $u$ is represented by $(a, x)$.  Similarly, one can fix
a representative $(b, y) \in {\bf R} \times Y_0$ of $v$, and express
$\Delta(u, v)$ as the infimum of $D(A(a, x), A(b, y))$ over all
representatives $(a, x) \in {\bf R} \times Y_0$ of $u$.  By fixing a
representative for $u$ or $v$ in ${\bf R} \times Y_0$, it is clear
that these infima are attained.  In particular, $\Delta(u, v) = 0$ if
and only if $u = v$, because of the analogous property of $D$ on ${\bf
R} \times Y_0$.

        In order to show that $\Delta(u, v)$ defines a metric on $Y$,
it remains to check that the triangle inequality holds, so that
\begin{equation}
\label{Delta(u, w) le Delta(u, v) + Delta(v, w)}
        \Delta(u, w) \le \Delta(u, v) + \Delta(v, w)
\end{equation}
for every $u, v, w \in Y$.  To do this, it suffices to verify that
\begin{equation}
\label{Delta(u, w) le D((a, x), (b, y)) + D((b', y'), (c, z))}
        \Delta(u, w) \le D((a, x), (b, y)) + D((b', y'), (c, z))
\end{equation}
for every $(a, x), (b, y), (b', y'), (c, z) \in {\bf R} \times Y_0$
such that $A(a, x) = u$, $A(b, y) = A(b', y') = v$, and $A(c, z) = w$.
The main point is to use translation-invariance of $D$ on ${\bf R}
\times Y_0$ to reduce to the case where $(b, y) = (b', y')$, as in the
previous paragraph.  In this case, we get that
\begin{equation}
\label{Delta(u, w) le ... le D((a, x), (b, y)) + D((b, y), (c, z))}
 \Delta(u, w) \le D((a, x), (c, z)) \le D((a, x), (b, y)) + D((b, y), (c, z))
\end{equation}
because of the triangle inequality for $D$ on ${\bf R} \times Y_0$, as
desired.  Observe also that
\begin{equation}
\label{Delta(A(a, x), A(b, y)) le D(A(a, x), A(b, y))}
        \Delta(A(a, x), A(b, y)) \le D(A(a, x), A(b, y))
\end{equation}
for every $(a, x), (b, y) \in {\bf R} \times Y_0$, by construction.

        If $x, y \in Y_0$, then we can identify $x$ and $y$ with
elements of $Y$ in the usual way, by putting $x_0 = y_0 = 0$ in ${\bf
R} / {\bf Z}$.  Let us check that
\begin{equation}
\label{Delta(x, y) = rho(x, y)}
        \Delta(x, y) = \rho(x, y).
\end{equation}
By definition, $\Delta(x, y)$ is equal to the infimum of $D(A(a, x),
A(b, y))$ over $a, b \in {\bf R}$ such that $q_0(a) = q_0(b) = 0$,
which is to say that $a, b \in {\bf Z}$.  If $a = b$, then
\begin{equation}
\label{D((a, x), (b, y)) = rho(x, y)}
        D((a, x), (b, y)) = \rho(x, y),
\end{equation}
and hence $\Delta(x, y) \le \rho(x, y)$.  Otherwise, if $a \ne b$,
then
\begin{equation}
\label{rho(x, y) le 1 le |a - b| le D((a, x), (b, y))}
        \rho(x, y) \le 1 \le |a - b| \le D((a, x), (b, y)),
\end{equation}
and (\ref{Delta(x, y) = rho(x, y)}) follows.

        Note that
\begin{equation}
\label{Delta(u, v) le 1}
        \Delta(u, v) \le 1
\end{equation}
for every $u, v \in Y$, since one can always choose $(a, x), (b, y) \in {\bf R}
\times Y_0$ such that $A(a, x) = u$, $A(b, y) = v$, and $|a - b| \le 1/2$.
Let us check that
\begin{equation}
\label{D((a, x), (b, y)) le 2 Delta(A(a, x), A(b, y))}
        D((a, x), (b, y)) \le 2 \, \Delta(A(a, x), A(b, y))
\end{equation}
for every $(a, x), (b, y) \in {\bf R} \times Y_0$ such that $|a - b|
\le 1/2$.  It suffices to show that
\begin{equation}
\label{D((a, x), (b, y)) le 2 D((a, x), (b', y'))}
        D((a, x), (b, y)) \le 2 \, D((a, x), (b', y'))
\end{equation}
whenever $(b', y') \in {\bf R} \times Y_0$ satisfies $A(b', y') = A(b,
y)$ and $(b', y') \ne (b, y)$.  In this case, $|b' - b| \ge 1$,
because $b' - b \in {\bf Z}$ and $b' \ne b$.  This implies that $|a -
b'| \ge 1/2$, and hence that
\begin{equation}
\label{D((a, x), (b, y)) le 1 le 2 |a - b'| le 2 D((a, x), (b, y))}
        D((a, x), (b, y)) \le 1 \le 2 \, |a - b'| \le 2 \, D((a, x), (b, y)),
\end{equation}
as desired.

\section{Haar measure on $Y$}
\label{haar measure on Y, 2}
\setcounter{equation}{0}

        Let $\mu$ be Haar measure on $Y$, normalized so that
\begin{equation}
\label{mu(Y) = 1}
        \mu(Y) = 1.
\end{equation}
Also let $p_n$ be the restriction of the coordinate projection $\pi_n$
in (\ref{pi_n(x) = x_n}) to $Y$ for each nonnegative integer $n$, so that
\begin{equation}
\label{p_n(x) = x_n}
        p_n(x) = x_n
\end{equation}
for every $x = \{x_l\}_{l = 0}^\infty \in Y$.  Thus $p_n$ is a
continuous homomorphism from $Y$ onto ${\bf R} / R_n \, {\bf Z}$ with
kernel equal to $Y_n$ for each $n \ge 0$.  If $E$ is a Borel
measurable subset of ${\bf R} / R_n \, {\bf Z}$, then $p_n^{-1}(E)$ is
a Borel measurable subset of $Y$, and
\begin{equation}
\label{mu(p_n^{-1}(E)) = |E|/R_n}
        \mu(p_n^{-1}(E)) = |E|/R_n.
\end{equation}
Here $|E|$ denotes the measure of $E$ as a subset of ${\bf R} / R_n \,
{\bf Z}$ that comes from Lebesgue measure on ${\bf R}$ in the obvious
way, so that $|{\bf R} / R_n \, {\bf Z}| = 1$.  If the $r_j$'s are
bounded, then one can check that $\mu$ is Ahlfors regular of dimension
$2$ on $Y$ with respect to the appropriate metrics discussed
previously.  Even if the $r_j$'s are not bounded, the subsets of $Y$
of the form $p_n^{-1}(E)$ with $E$ a Borel set in ${\bf R} / R_n \,
{\bf Z}$ define a nice filtration on $Y$, with the corresponding
martingales on $Y$.

\part{Another perspective}
\label{another perspective}

\section{Ultrametrics}
\label{ultrametrics}
\setcounter{equation}{0}

        A metric $d(x, y)$ on a set $M$ is said to be an \emph{ultrametric} if
\begin{equation}
\label{d(x, z) le max(d(x, y), d(y, z))}
        d(x, z) \le \max(d(x, y), d(y, z))
\end{equation}
for every $x, y, z \in M$, which is stronger than the usual triangle
inequality.  The discrete metric on any set is an ultrametric, for
instance.  

        As another class of examples, let $X_1, X_2, X_3, \ldots$ be a 
sequence of nonempty sets, and let $X = \prod_{j = 1}^\infty X_j$ be
their Cartesian product, consisting of all sequences $x = \{x_j\}_{j =
  1}^\infty$ such that $x_j \in X_j$ for each $j$.  If $x, y \in X$,
then let $n(x, y)$ be the largest nonnegative integer such that $x_j =
y_j$ for each $j \le n(x, y)$, with $n(x, y) = +\infty$ when $x = y$.
Equivalently, $n(x, y) + 1$ is the smallest positive integer $j$ such
that $x_j \ne y_j$ when $x \ne y$.  Of course,
\begin{equation}
\label{n(x, y) = n(y, x)}
        n(x, y) = n(y, x)
\end{equation}
for every $x, y \in X$, and it is easy to see that
\begin{equation}
\label{n(x, z) ge min(n(x, y), n(y, z))}
        n(x, z) \ge \min(n(x, y), n(y, z))
\end{equation}
for every $x, y, z \in X$.  Let $t = \{t_j\}_{j = 0}^\infty$ be a
decreasing sequence of positive real numbers that converges to $0$.
Put
\begin{equation}
\label{d(x, y) = t_{n(x, y)}}
        d(x, y) = t_{n(x, y)}
\end{equation}
for each $x, y \in X$ with $x \ne y$, and $d(x, y) = 0$ when $x = y$,
which corresponds to taking $t_\infty = 0$.  It is easy to see that
$d(x, y)$ is an ultrametric on $X$, using (\ref{n(x, y) = n(y, x)})
and (\ref{n(x, z) ge min(n(x, y), n(y, z))}).  The topology on $X$
determined by this ulrametric is the same as the product topology
corresponding to the discrete topology on $X_j$ for each $j$.  In
particular, $X$ is compact with respect to this topology when $X_j$ is
a finite set for each $j$.

        If $d(x, y)$ is a metric on any set $M$, then the open ball
in $M$ centered at a point $x \in M$ and with radius $r > 0$ is defined
as usual by
\begin{equation}
\label{B(x, r) = {y in M : d(x, y) < r}}
        B(x, r) = \{y \in M : d(x, y) < r\}.
\end{equation}
If $y \in B(x, r)$, then $t = r - d(x, y) > 0$, and one can check that
\begin{equation}
\label{B(y, t) subseteq B(x, r)}
        B(y, t) \subseteq B(x, r),
\end{equation}
using the triangle inequality.  However, if $d(\cdot, \cdot)$ is an
ultrametric on $M$, then
\begin{equation}
\label{B(y, r) subseteq B(x, r)}
        B(y, r) \subseteq B(x, r)
\end{equation}
for every $y \in B(x, r)$.  This implies that
\begin{equation}
\label{B(x, r) = B(y, r)}
        B(x, r) = B(y, r)
\end{equation}
when $d(x, y) < r$, since the same argument can be applied with the
roles of $x$ and $y$ reversed.

        Similarly, the closed ball in $M$ 
\begin{equation}
\label{overline{B}(x, r) = {y in M : d(x, y) le r}}
        \overline{B}(x, r) = \{y \in M : d(x, y) \le r\},
\end{equation}
with center $x \in M$ and radius $r \ge 0$ is a closed set in $M$ with
respect to the topology determined by $d(\cdot, \cdot)$, for any
metric $d(\cdot, \cdot)$ on $M$.  If $d(\cdot, \cdot)$ is an
ultrametric on $M$, then
\begin{equation}
\label{overline{B}(y, r) subseteq overline{B}(x, r)}
        \overline{B}(y, r) \subseteq \overline{B}(x, r)
\end{equation}
for every $y \in \overline{B}(x, r)$, which implies that
$\overline{B}(x, r)$ is also an open set in $M$.  As before, one can
apply this with the roles of $x$ and $y$ reversed, to get that
\begin{equation}
\label{overline{B}(x, r) = overline{B}(y, r)}
        \overline{B}(x, r) = \overline{B}(y, r)
\end{equation}
when $d(x, y) \le r$.

        If $d(\cdot, \cdot)$ is an ultrametric on $M$, $x, y, z \in M$,
and $d(y, z) \le d(x, y)$, then
\begin{equation}
\label{d(x, z) le max(d(x, y), d(y, z)) = d(x, y)}
        d(x, z) \le \max(d(x, y), d(y, z)) = d(x, y).
\end{equation}
If $d(y, z) < d(x, y)$, then
\begin{equation}
\label{d(x, y) le max(d(x, z), d(y, z))}
        d(x, y) \le \max(d(x, z), d(y, z))
\end{equation}
implies that
\begin{equation}
\label{d(x, y) le d(x, z)}
        d(x, y) \le d(x, z).
\end{equation}
Hence
\begin{equation}
\label{d(x, y) = d(x, z)}
        d(x, y) = d(x, z)
\end{equation}
when $d(y, z) < d(x, y)$, by (\ref{d(x, z) le max(d(x, y), d(y, z)) =
  d(x, y)}).

        Let $d(\cdot, \cdot)$ be a metric on a set $M$ again, and put
\begin{equation}
\label{V(x, r) = {y in M : d(x, y) > r}}
        V(x, r) = \{y \in M : d(x, y) > r\}
\end{equation}
for each $x \in M$ and $r \ge 0$, which is the same as $M \setminus
\overline{B}(x, r)$.  If $y \in V(x, r)$, then $t = d(x, y) - r > 0$,
and one can check that
\begin{equation}
\label{B(y, t) subseteq V(x, r)}
        B(y, t) \subseteq V(x, r),
\end{equation}
using the triangle inequality.  If $d(\cdot, \cdot)$ is an ultrametric
on $M$, then
\begin{equation}
\label{B(y, d(x, y)) subseteq V(x, r)}
        B(y, d(x, y)) \subseteq V(x, r)
\end{equation}
for every $y \in V(x, r)$, by (\ref{d(x, y) = d(x, z)}).

        Similarly, put
\begin{equation}
\label{W(x, r) = {y in M : d(x, y) ge r}}
        W(x, r) = \{y \in M : d(x, y) \ge r\}
\end{equation}
for each $x \in M$ and $r > 0$, which is the same as $M \backslash
B(x, r)$.  If $d(\cdot, \cdot)$ is an ultrametric on $M$, then
\begin{equation}
\label{B(y, r) subseteq B(y, d(x, y)) subseteq W(x, r)}
        B(y, r) \subseteq B(y, d(x, y)) \subseteq W(x, r)
\end{equation}
for every $y \in W(x, r)$, by (\ref{d(x, y) = d(x, z)}).  In
particular, this implies that $W(x, r)$ is an open set, so that $B(x,
r)$ is a closed set in $M$.

\section{$r$-Adic absolute values}
\label{r-adic absolute values}
\setcounter{equation}{0}

        Let $r = \{r_j\}_{j = 1}^\infty$ be a sequence of positive integers,
with $r_j \ge 2$ for each $j$.  Put
\begin{equation}
\label{R_l = prod_{j = 1}^l r_j, 2}
        R_l = \prod_{j = 1}^l r_j
\end{equation}
for each positive integer $l$ and $R_0 = 1$, and let $t = \{t_l\}_{l =
  0}^\infty$ be a decreasing sequence of real numbers that converges
to $0$.  One can take $t_l = 1/R_l$ for each $l$, for instance, and in
any case one might at least take $t_0 = 1$.  If $r_j = (r_1)^j$ for
each $j \ge 1$, then $R_l = (r_1)^l$ for each $l \ge 0$, and $t_l =
1/R_l = (r_1)^{-l}$ is especially nice.  In particular, if $p$ is a
prime number, $r_j = p$ for every $j \ge 1$, and $t_l = p^{-l}$ for
each $l \ge 0$, then this reduces to the usual situation for $p$-adic
numbers.

        Let $a$ be an integer, and let $n(a)$ be the largest nonnegative
integer $l$ such that $a$ is an integer multiple of $R_l$, with 
$n(0) = +\infty$.  Thus
\begin{equation}
\label{n(-a) = n(a)}
        n(-a) = n(a)
\end{equation}
for each integer $a$, and it is easy to see that
\begin{equation}
\label{n(a + b) ge min(n(a), n(b))}
        n(a + b) \ge \min(n(a), n(b))
\end{equation}
and
\begin{equation}
\label{n(a b) ge max(n(a), n(b))}
        n(a \, b) \ge \max(n(a), n(b))
\end{equation}
for any two integers $a$, $b$.  If $r_j = r_1$ for each $j \ge 1$,
then
\begin{equation}
\label{n(a b) ge n(a) + n(b)}
        n(a \, b) \ge n(a) + n(b)
\end{equation}
for all integers $a$, $b$.  If $r_j = p^j$ for some prime number $p$
and every $j \ge 1$, then
\begin{equation}
\label{n(a b) = n(a) + n(b)}
        n(a \, b) = n(a) + n(b)
\end{equation}
for each $a$, $b$.

        The $r$-adic absolute value of an integer $x$ is defined by
\begin{equation}
\label{|a|_r = t_{n(a)}}
        |a|_r = t_{n(a)}
\end{equation}
when $a \ne 0$ and $|0|_r = 0$, which corresponds to (\ref{|a|_r =
  t_{n(a)}}) with $t_\infty = 0$.  Note that
\begin{equation}
\label{|a|_r le t_0}
        |a|_r \le t_0
\end{equation}
for every integer $x$, and that
\begin{equation}
\label{|-a|_r = |a|_r}
        |-a|_r = |a|_r
\end{equation}
by (\ref{n(-a) = n(a)}).  Similarly,
\begin{equation}
\label{|a + b|_r le max(|a|_r, |b|_r), 2}
        |a + b|_r \le \max(|a|_r, |b|_r)
\end{equation}
and
\begin{equation}
\label{|a b|_r le min(|a|_r, |b|_r)}
        |a \, b|_r \le \min(|a|_r, |b|_r)
\end{equation}
for every $a$ and $b$, by (\ref{n(a + b) ge min(n(a), n(b))}) and
(\ref{n(a b) ge max(n(a), n(b))}).  Suppose that $r_j = r_1$ for each
$j \ge 1$ and $t = \{t_l\}_{l = 0}^\infty$ is submultiplicative, in
the sense that
\begin{equation}
\label{t_{k + l} le t_k t_l}
        t_{k + l} \le t_k \, t_l
\end{equation}
for every $k, l \ge 0$.  Under these conditions, we get that
\begin{equation}
\label{|a b|_r le |a|_r |b|_r, 2}
        |a \, b|_r \le |a|_r \, |b|_r
\end{equation}
for every $a$ and $b$, by (\ref{n(a b) ge n(a) + n(b)}).  If $r_j = p$
for some prime number $p$ and every $j$, and if $t_l = (t_1)^l$ for
each $l \ge 0$, then
\begin{equation}
\label{|a b|_r = |a|_r |b|_r}
        |a \, b|_r = |a|_r \, |b|_r
\end{equation}
for every $a$ and $b$, by (\ref{n(a b) = n(a) + n(b)}).  In
particular, $|a|_r$ reduces to the usual $p$-adic absolute value
$|a|_p$ of $a$ when $p$ is a prime number, $r_j = p$ for each $j \ge
1$, and $t_l = p^{-l}$ for every $l \ge 0$.

        The $r$-adic metric on the set ${\bf Z}$ of integers is defined by
\begin{equation}
\label{d_r(a, b) = |a - b|_r}
        d_r(a, b) = |a - b|_r.
\end{equation}
It is easy to see that this defines an ultrametric on ${\bf Z}$, using
(\ref{|-a|_r = |a|_r}) and (\ref{|a + b|_r le max(|a|_r, |b|_r), 2}).
Note that the topology determined on ${\bf Z}$ by (\ref{d_r(a, b) = |a
  - b|_r}) depends only on $r$, and not on $t$.  Thus we shall
sometimes refer to this topology simply as the $r$-adic topology on
${\bf Z}$.  If $r_j = p$ for some prime number $p$ and every $j \ge
1$, and if $t_l = p^{-l}$ for each $l \ge 0$, then the $r$-adic metric
on ${\bf Z}$ reduces to the usual $p$-adic metric.

\section{Coherent sequences}
\label{coherent sequences}
\setcounter{equation}{0}

        Let us continue with the same notation and hypotheses as before,
and put
\begin{equation}
\label{X_l = {bf Z} / R_l {bf Z}}
        X_l = {\bf Z} / R_l \, {\bf Z}
\end{equation}
for each positive integer $l$.  Thus $X_l$ is a commutative ring with
$R_l$ elements for each $l$.  Consider the Cartesian product
\begin{equation}
\label{X = prod_{l = 1}^infty X_l}
        X = \prod_{l = 1}^\infty X_l,
\end{equation}
which is also a commutative ring with respect to coordinatewise
addition and multiplication.  Note that $X$ is a compact Hausdorff
topological space, with respect to the product topology corresponding
to the discrete topology on $X_l$ for each $l$.  It is easy to see that
addition and multiplication on $X$ are continuous with respect to the
product topology, so that $X$ is a topological ring.

        Let $q_l$ be the usual quotient homomorphism from ${\bf Z}$
onto ${\bf Z} / R_l \, {\bf Z}$ for each $l \ge 1$.  If we put
\begin{equation}
\label{q(a) = {q_l(a)}_{l = 1}^infty}
        q(a) = \{q_l(a)\}_{l = 1}^\infty
\end{equation}
for each $a \in {\bf Z}$, then $q$ defines a ring homomorphism from
${\bf Z}$ into $X$.  It is easy to see that the kernel of this
homomorphism is trivial, since $R_l \to \infty$ as $l \to \infty$.

        Observe that
\begin{equation}
\label{n(q(a), q(b)) = n(a - b)}
        n(q(a), q(b)) = n(a - b)
\end{equation}
for every $a, b \in {\bf Z}$, where $n(q(a), q(b))$ is as in Section
\ref{ultrametrics}, and $n(a - b)$ is as in Section \ref{r-adic
  absolute values}.  This implies that
\begin{equation}
\label{d(q(a), q(b)) = d_r(a, b)}
        d(q(a), q(b)) = d_r(a, b)
\end{equation}
for every $a, b \in {\bf Z}$, where $d(\cdot, \cdot)$ is defined on
$X$ as in (\ref{d(x, y) = t_{n(x, y)}}), and $d_r(a, b)$ is the
$r$-adic metric on ${\bf Z}$, as in (\ref{d_r(a, b) = |a - b|_r}).  Of
course, it is important to use the same sequence $t = \{t_l\}_{l =
  0}^\infty$ in both cases.  In particular, $q$ is a homeomorphism
from ${\bf Z}$ onto its image in $X$, where ${\bf Z}$ is equipped with
the topology determined by the $r$-adic metric, and $X$ is equipped
with the product topology mentioned earlier.

        There is a natural ring homomorphism from $X_{l + 1}$ onto $X_l$ 
for each $l$, because $R_{l + 1} = r_{l + 1} \, R_l$, and hence $R_{l
  + 1} \, {\bf Z} \subseteq R_l \, {\bf Z}$.  Equivalently, this
homomorphism maps $q_{l + 1}(a) \in {\bf Z} / R_{l + 1} \, {\bf Z}$ to
$q_l(a) \in {\bf Z} / R_l \, {\bf Z}$ for each $a \in {\bf Z}$.  An
element $x = \{x_l\}_{l = 1}^\infty$ of $X$ is said to be a
\emph{coherent sequence} if $x_l$ is the image of $x_{l + 1}$ under
this homomorphism from $X_{l + 1}$ onto $X_l$ for each $l$.  Thus
$q(a)$ is a coherent sequence for each $a \in {\bf Z}$.  It is easy to
see that the set of coherent sequences forms a sub-ring of $X$, and
also a closed set in $X$ with respect to the product topology.

        In fact, the set of coherent sequences in $X$ is the same as
the closure of $q({\bf Z})$ in $X$ with respect to the product
topology.  To see this, it suffices to check that every coherent
sequence $x \in X$ can be approximated by elements of $q({\bf Z})$
with respect to the product topology on $X$.  Of course, for any $x
\in X$ and positive integer $k$, there is an $a \in {\bf Z}$ such that
$q_k(a) = x_k$.  If $x$ is a coherent sequence, then it follows that
$q_j(a) = x_j$ for every $j \le k$, which implies that $x$ can be
approximated by elements of $q({\bf Z})$ with respect to the product
topology, as desired, by taking $k \to \infty$.

        It is easy to check directly that $X$ is complete with respect
to the metric $d(x, y)$ in (\ref{d(x, y) = t_{n(x, y)}}), and hence that
the set of coherent sequences is also complete with respect to this
metric.  Because $q$ is an isometric embedding of ${\bf Z}$ into $X$,
as in (\ref{d(q(a), q(b)) = d_r(a, b)}), the completion of ${\bf Z}$
with respect to the $r$-adic metric can be identified with the set of
coherent sequences in $X$.  Let us refer to this completion as the
ring ${\bf Z}_r$ of $r$-adic integers.  Thus ${\bf Z}_r$ is a compact
commutative topological ring which contains ${\bf Z}$ as a dense sub-ring,
and the $r$-adic metric extends to ${\bf Z}_r$ in a natural way.
The $r$-adic absolute value function also extends to ${\bf Z}_r$
in a natural way, and satisfies properties like those in the previous
section.

        Note that the identification of $r$-adic integers with coherent
sequences does not depend on the choice of sequence $t = \{t_l\}_{l =
  0}^\infty$ in the definition of the $r$-adic absolute value function
and metric.  Different choices of $t$ lead to topologically-equivalent
translation-invariant $r$-adic metrics on ${\bf Z}$ anyway, which have
the same Cauchy sequences, and isomorphic completions.  If $p$ is a
prime number, $r_j = p$ for each $j \ge 1$, and $t_l = p^{-l}$ for
each $l \ge 0$, then this description of ${\bf Z}_r$ is equivalent to
the usual ring ${\bf Z}_p$ of $p$-adic integers.

\section{Topological equivalence}
\label{topological equivalence}
\setcounter{equation}{0}

        Let $r = \{r_j\}_{j = 1}^\infty$ and $r' = \{r_j'\}_{j = 1}^\infty$
be sequences of integers, with $r_j, r_j' \ge 2$ for each $j$.  As
before, put
\begin{equation}
\label{R_l = prod_{j = 1}^l r_j, R_l' = prod_{j = 1}^l r_j'}
        R_l = \prod_{j = 1}^l r_j, \quad R_l' = \prod_{j = 1}^l r_j'
\end{equation}
when $l \ge 1$, and $R_0 = R_0' = 1$.  If for each $k \ge 1$ there is
an $l \ge 1$ such that $R_l'$ is an integer multiple of $R_k$, then
put $r \prec r'$.  If $r \prec r'$ and $r' \prec r$, then put $r \sim
r'$.  The relation $r \prec r'$ is clearly reflexive and transitive,
which implies that $r \sim r'$ is an equivalence relation.

        Note that $\{R_l\}_{l = 1}^\infty$ automatically converges to $0$
with respect to the $r$-adic topology on ${\bf Z}$, and similarly that
$\{R_l'\}_{l = 1}^\infty$ automatically converges to $0$ with respect
to the $r'$-adic topology on ${\bf Z}$.  It is easy to see that
$\{R_l'\}_{l = 1}^\infty$ converges to $0$ with respect to the
$r$-adic topology on ${\bf Z}$ if and only if $r \prec r'$.  If $r
\prec r'$, then the $r$-adic topology on ${\bf Z}$ is weaker than the
$r'$-adic topology on ${\bf Z}$.  Conversely, if the $r$-adic topology
on ${\bf Z}$ is weaker than the $r'$-adic topology, then $\{R_l'\}_{l
  = 1}^\infty$ converges to $0$ with respect to the $r$-adic topology,
and hence $r \prec r'$.  In this case, every coherent sequence with
respect to $r'$ determines a coherent sequence with respect to $r$,
which leads to a continuous ring homomorphism from ${\bf Z}_{r'}$ into
${\bf Z}_r$.  One can also check that this mapping is surjective,
because ${\bf Z}$ is dense in ${\bf Z}_r$, and ${\bf Z}_{r'}$ is
compact.  It follows that $r \sim r'$ if and only if the $r$-adic and
$r'$-adic topologies on ${\bf Z}$ are the same, in which event ${\bf
  Z}_r$ and ${\bf Z}_{r'}$ are isomorphic as topological rings.

        Let $p$ be a prime number, and let $c_{r, l}(p)$ be the number
of factors of $p$ in $R_l$ for each positive integer $l$.  Thus $c_{r,
  l}(p) \le c_{r, l + 1}(p)$ for each $l$, and we put
\begin{equation}
\label{c_r(p) = sup_{l ge 1} c_{r, l}(p)}
        c_r(p) = \sup_{l \ge 1} c_{r, l}(p),
\end{equation}
which is either a nonnegative integer or $+\infty$.  If $c_{r'}(p)$ is
defined in the same way, then $r \prec r'$ if and only if $c_r(p) \le
c_{r'}(p)$ for every prime number $p$, and hence $r \sim r'$ if and
only if $c_r(p) = c_{r'}(p)$ for every prime number $p$.  Because the
number of factors in $R_l$ is strictly increasing, either $c_r(p) =
+\infty$ for some $p$, or $c_r(p) > 0$ for infinitely many $p$.
Conversely, if $c(p)$ is any function on the set of prime numbers such
that $c(p)$ is a nonnegative integer or $+\infty$ for each $p$, and
either $c(p) = +\infty$ for some $p$ or $c(p) > 0$ for infinitely many
$p$, then $c(p) = c_r(p)$ for some $r$ as before.

        Consider the Cartesian product
\begin{equation}
\label{(prod_{0 < c_r(p) < infty} ... times (prod_{c_r(p) = +infty} ...)}
        \Big(\prod_{0 < c_r(p) < \infty} ({\bf Z} / p^{c_r(p)} \, {\bf Z})\Big)
         \times \Big(\prod_{c_r(p) = +\infty} {\bf Z}_p\Big),
\end{equation}
where more precisely one takes the product over the prime numbers $p$
such that $0 < c_r(p) < \infty$ and $c_r(p) = +\infty$, respectively.
This is a commutative ring with respect to coordinatewise addition and
multiplication.  This is also a compact topological ring with respect
to the product topology corresponding to the discrete topology on
${\bf Z} / p^{c_r(p)} \, {\bf Z}_p$ when $0 < c_r(p) < \infty$, and
the $p$-adic topology on ${\bf Z}_p$ when $c_r(p) = \infty$.  As
before, there is a natural continuous ring homomorphism from ${\bf
  Z}_r$ onto each of the factors, which one can get using coherent
sequences.  This leads to a continuous ring homomorphism from ${\bf
  Z}_r$ into (\ref{(prod_{0 < c_r(p) < infty} ... times (prod_{c_r(p)
    = +infty} ...)}).  The image of ${\bf Z}$ is dense in
(\ref{(prod_{0 < c_r(p) < infty} ... times (prod_{c_r(p) = +infty}
  ...)}) with respect to the product topology, by the chinese
remainder theorem.  This implies that ${\bf Z}_r$ maps onto
(\ref{(prod_{0 < c_r(p) < infty} ... times (prod_{c_r(p) = +infty}
  ...)}), because ${\bf Z}_r$ is compact.  The kernel of this
homomorphism from ${\bf Z}_r$ onto (\ref{(prod_{0 < c_r(p) < infty}
  ... times (prod_{c_r(p) = +infty} ...)}) is trivial, and this
isomorphism is a homeomorphism as well.

\section{Comparisons}
\label{comparisons}
\setcounter{equation}{0}

        Let $X_1, X_2, X_3, \ldots$ be a sequence of nonempty finite sets,
each of which has at least two elements, and let $X = \prod_{j = 1}^\infty X_j$ 
be their Cartesian product, as in Section \ref{ultrametrics}.  Also
let $\{t_j\}_{j = 0}^\infty$ be a strictly decreasing sequence of
positive real numbers that converges to $0$, which leads to an
ultrametric $d(x, y)$ on $X$ as in (\ref{d(x, y) = t_{n(x, y)}}).  As
before, the topology on $X$ determined by $d(x, y)$ is the same as the
product topology corresponding to the discrete topology on each $X_j$.
In particular, $X$ is a compact Hausdorff space, and more precisely a
topological Cantor set.

        Suppose that $\{\widetilde{t}_j\}_{j = 1}^\infty$ is another
strictly decreasing sequence of positive real numbers that converges
to $0$, and let $\widetilde{d}(x, y)$ be the corresponding ultrametric
on $X$, as in (\ref{d(x, y) = t_{n(x, y)}}).  Thus $d(x, y)$ and 
$\widetilde{d}(x, y)$ determine the same topology on $X$, and in fact
they determine the same collections of open and closed balls in $X$.
This is a very strong geometric property, and indeed this collection
of balls has a lot of interesting structure.  The nesting of these balls
leads to very simple covering lemmas, which lead in turn to maximal
function estimates.

        Now let $\pi$ be a one-to-one mapping of the set ${\bf Z}_+$
of positive integers onto itself.  This can be used to rearrange the
initial sequence of $X_j$'s to get a new sequence $X_{\pi(1)}, X_{\pi(2)}, 
X_{\pi(3)}, \ldots$ of sets, and thus a new Cartesian product 
$X^\pi = \prod_{j = 1}^\infty X_{\pi(j)}$.  Of course, there is also a
natural one-to-one correspondence between $X$ and $X^\pi$, which
sends $x = \{x_j\}_{j = 1}^\infty \in X$ to $\{x_{\pi(j)}\}_{j = 1}^\infty
\in X^\pi$.  This mapping is a homeomorphism between $X$ and $X^\pi$,
with respect to the product topologies on $X$ and $X^\pi$ corresponding
to the discrete topologies on the $X_j$'s.  However, this type of mapping
can still change the geometric structures being considered in significant
ways.

        Let $r = \{r_j\}_{j = 1}^\infty$ be a sequence of integers
with $r_j \ge 2$ for each $j$ again, and let $t = \{t_l\}_{l = 0}^\infty$
be a strictly decreasing sequence of positive real numbers that converges
to $0$.  As before, this leads to an $r$-adic absolute value and metric
on ${\bf Z}$, and on the corresponding completion ${\bf Z}_r$.  If
$\{\widetilde{t}_l\}_{l = 0}^\infty$ is another strictly decreasing
sequence of positive real numbers converging to $0$, then one gets
another $r$-adic absolute value function and metric, but the same topology
and completion ${\bf Z}_r$.  One also gets the same collections of
open and closed balls in ${\bf Z}$ and ${\bf Z}_r$, as in the previous
situation.  If $r' = \{r_j'\}_{j = 1}^\infty$ is another sequence with
$r \sim r'$, as in the previous section, then one gets the same topology
on ${\bf Z}$ and an isomorphic completion ${\bf Z}_{r'}$, but the corresponding
geometry can be affected significantly.

\section{Solenoids}
\label{solenoids}
\setcounter{equation}{0}

        Let $r = \{r_j\}_{j = 1}^\infty$ be a sequence of integers with
$r_j \ge 2$ for each $j$, and let $R_l$ be as in 
(\ref{R_l = prod_{j = 1}^l r_j, 2}).  This leads to the ring ${\bf Z}_r$
of $r$-adic integers, as before, which contains ${\bf Z}$ as a dense
sub-ring.  We can also consider ${\bf Z}$ as a discrete subgroup of
the real line ${\bf R}$ with respect to addition, with the
corresponding quotient ${\bf R} / {\bf Z}$ as a compact commutative
topological group.  Of course, ${\bf R} \times {\bf Z}_r$ is a locally
compact commutative topological group with respect to coordinatewise
addition and the product topology, using the standard topology on
${\bf R}$.  The group ${\bf Z}$ of integers with respect to addition
is a subgroup of both ${\bf R}$ and ${\bf Z}_r$, and hence
\begin{equation}
\label{A = {(a, a) : a in {bf Z}}}
        A = \{(a, a) : a \in {\bf Z}\}
\end{equation}
is a subgroup of ${\bf R} \times {\bf Z}_r$.  More precisely, $A$ is a
closed subgroup of ${\bf R} \times {\bf Z}_r$, because ${\bf Z}$ is a
closed subgroup of ${\bf R}$, so that the quotient group
\begin{equation}
\label{({bf R} times {bf Z}_r) / A}
        ({\bf R} \times {\bf Z}_r) / A
\end{equation}
is also a topological group with respect to the quotient topology.
The canonical quotient mapping from ${\bf R} \times {\bf Z}_r$ onto
(\ref{({bf R} times {bf Z}_r) / A}) is a local homeomorphism, and it
is easy to see that (\ref{({bf R} times {bf Z}_r) / A}) is compact,
because ${\bf Z}_r$ is compact.  One can also check that the image of
${\bf R} \times \{0\}$ in (\ref{({bf R} times {bf Z}_r) / A}) under
the canonical quotient mapping from ${\bf R} \times {\bf Z}_r$ onto
(\ref{({bf R} times {bf Z}_r) / A}) is dense in (\ref{({bf R} times
  {bf Z}_r) / A}), because ${\bf Z}$ is dense in ${\bf Z}_r$.  This
implies that (\ref{({bf R} times {bf Z}_r) / A}) is connected, because
${\bf R}$ is connected, and the closure of a connected set is
connected.

        There is also a nice description of (\ref{({bf R} times {bf Z}_r) / A})
in terms of coherent sequences.  Put
\begin{equation}
\label{Y_l = {bf R} / R_l {bf Z}}
        Y_l = {\bf R} / R_l \, {\bf Z}
\end{equation}
for each $l \ge 0$, considered as a compact commutative topological
group.  Thus
\begin{equation}
        Y = \prod_{l = 0}^\infty Y_l,
\end{equation}
is also a compact commutative topological group, with respect to
coordinatewise addition and the product topology.  There is a natural
group homomorphism from $Y_{l + 1}$ onto $Y_l$ for each $l \ge 0$,
because $R_{l + 1} \, {\bf Z} \subseteq R_l \, {\bf Z}$.  This
homomorphism is also continuous, and in fact a local homeomorphism.
An element $y = \{y_l\}_{l = 0}^\infty$ of $Y$ is said to be a
coherent sequence if $y_l$ is the image of $y_{l + 1}$ under this
homomorphism for each $l \ge 0$.  The set of coherent sequences in $Y$
is a closed subgroup of $Y$, and there is a natural isomorphism
between ${\bf Z}_r$ and the coherent sequences $y = \{y_l\}_{l =
  0}^\infty$ in $Y$ such that $y_0 = 0$ in $Y_0 = {\bf R} / {\bf Z}$.
There is a natural homomorphism from ${\bf R}$ into $Y$, whose $l$th
coordinate is the canonical quotient homomorphism from ${\bf R}$ onto
$Y_l$ for each $l \ge 0$, and it is easy to see that this homomorphism
sends real numbers to coherent sequences in $Y$.  This leads to a
homomorphism from ${\bf R} \times {\bf Z}_r$ into the group of
coherent sequences in $Y$, which adds the images of elements of ${\bf
  R}$ and ${\bf Z}_r$ in $Y$.  One can check that the kernel of this
homomorphism is equal to (\ref{A = {(a, a) : a in {bf Z}}}), which
leads to an isomorphism between (\ref{({bf R} times {bf Z}_r) / A})
and the group of coherent sequences in $Y$ as topological groups.

        If $t = \{t_l\}_{l = 0}^\infty$ is a decreasing sequence of
positive real numbers that converges to $0$, then we get a
corresponding $r$-adic metric on ${\bf Z}_r$, which is invariant under
translations.  The standard metric on ${\bf R}$ is invariant under
translations too, and one can combine the two metrics to get a
translation-invariant metric on ${\bf R} \times {\bf Z}_r$.  Using
this, one can get a translation-invariant quotient metric on
(\ref{({bf R} times {bf Z}_r) / A}).

        Suppose that $r' = \{r_j'\}_{j = 1}^\infty$ is another sequence
of integers with $r_j' \ge 2$ for each $j$ such that $r \sim r'$, in
the sense of Section \ref{topological equivalence}.  It is easy to see
that (\ref{({bf R} times {bf Z}_r) / A}) is isomorphic as a
topological group to its analogue with $r'$ instead of $r$, because of
the isomorphism between ${\bf Z}_r$ and ${\bf Z}_{r'}$, which is the
identity mapping on their common subgroup ${\bf Z}$.  One can also
look at this in terms of coherent sequences.  However, this type of
isomorphism may be rather complicated with respect to the corresponding
geometries.

\section{Filtrations}
\label{filtrations}
\setcounter{equation}{0}

        Remember that a filtration on a probability space
$(X, \mathcal{A}, \mu)$ is an increasing sequence
$\mathcal{B}_1 \subseteq \mathcal{B}_2 \subseteq \mathcal{B}_3 \cdots$
of $\sigma$-subalgebras of $\mathcal{A}$.  This can be used to define
conditional expectation operators, martingales, and so on.  As a basic
class of examples, suppose that $(X_j, \mathcal{A}_j, \mu_j)$ is a
probability space for each positive integer $j$, and let
\begin{equation}
\label{X= prod_{j  = 1}^infty X_j}
        X= \prod_{j  = 1}^\infty X_j
\end{equation}
be their product, with the corresponding product $\sigma$-algebra
$\mathcal{A}$ and probability measure $\mu$.  Let $\mathcal{B}_l$ be
the collection of subsets of $\prod_{j = 1}^\infty X_j$ that
correspond to a product of a measurable subset of $\prod_{j = 1}^l
X_j$ with $\prod_{j = l + 1}^\infty X_j$, for each positive integer
$l$.  It is easy to see that this defines a filtration on $X$.  If
$\pi$ is a one-to-one mapping from ${\bf Z}_+$ onto itself, then one
can also consider the product $X^\pi = \prod_{j = 1}^\infty
X_{\pi(j)}$ as a probability space, with the analogous filtration.
The mapping from $x = \{x_j\}_{j = 1}^\infty \in X$ to
$\{x_{\pi(j)}\}_{j = 1}^\infty \in X^\pi$ defines an isomorphism
between $X$ and $X^\pi$ as probability spaces, but this isomorphism
may be rather complicated in terms of the corresponding filtrations.
In particular, one might take $X_j$ to be a finite set with at least
two elements for each $j$, where every subset of $X_j$ is measurable,
and where $\mu_j$ assigns equal weight to each element of $X_j$.  The
resulting filtration on (\ref{X= prod_{j = 1}^infty X_j}) is closely
related to the type of ultrametrics on $X$ discussed earlier.

        Now let $r = \{r_j\}_{j = 1}^\infty$ be a sequence of integers
with $r_j \ge 2$ for each $j$, and let $R_l$ be as in 
(\ref{R_l = prod_{j = 1}^l r_j, 2}).  Because ${\bf Z}_r$ is a compact
commutative topological group, there is a natural translation-invariant
Borel probability measure on ${\bf Z}_r$, given by Haar measure.
If $r' = \{r_j'\}_{j = 1}^\infty$ is another sequence of integers
with $r_j' \ge 2$ for each $j$ and $r \sim r'$, then ${\bf Z}_r$
is isomorphic to ${\bf Z}_{r'}$ as a topological ring and hence as a
topological group, and Haar measure on ${\bf Z}_r$ corresponds to
Haar measure on ${\bf Z}_{r'}$ under this isomorphism.

        Of course, there is a natural homomorphism from ${\bf Z}$
onto ${\bf Z} / R_l \, {\bf Z}$ for each $l$.  This extends to a
continuous homomorphism from ${\bf Z}_r$ onto ${\bf Z} / R_l \, {\bf Z}$
for each $l$, basically by construction.  Let $\mathcal{B}_l$
be the collection of subsets of ${\bf Z}_r$ that can be expressed as
the inverse image of a subset of ${\bf Z} / R_l \, {\bf Z}$ under the
homomorphism from ${\bf Z}_r$ onto ${\bf Z} / R_l \, {\bf Z}$ just
mentioned.  This is a $\sigma$-subalgebra of the Borel sets in ${\bf Z}_r$
for each $l$, which defines a filtration on ${\bf Z}_r$.  This filtration
is closely related to the $r$-adic geometry on ${\bf Z}_r$.

        Similarly, (\ref{({bf R} times {bf Z}_r) / A}) is a compact
commutative topological group, and thus has a natural
translation-invariant Borel probability measure, given by Haar
measure.  Remember that there is a natural continuous homomorphism
from (\ref{({bf R} times {bf Z}_r) / A}) onto (\ref{Y_l = {bf R} / R_l
  {bf Z}}) for each nonnegative integer $l$, as in the previous
section.  Let $\mathcal{C}_l$ be the collection of subsets of
(\ref{({bf R} times {bf Z}_r) / A}) which can be expressed as the
inverse image of a Borel subset of (\ref{Y_l = {bf R} / R_l {bf Z}})
under the homomorphism from (\ref{({bf R} times {bf Z}_r) / A}) onto
(\ref{Y_l = {bf R} / R_l {bf Z}}) just mentioned.  This is a
$\sigma$-subalgebra of the Borel sets in (\ref{({bf R} times {bf Z}_r)
  / A}) for each $l$, which defines a filtration on (\ref{({bf R}
  times {bf Z}_r) / A}).  In effect, this filtration is part of the
geometry of (\ref{({bf R} times {bf Z}_r) / A}), which is specified by
$r$.

\section{Directed systems}
\label{directed systems}
\setcounter{equation}{0}

        Let $(A, \preceq)$ be a partially-ordered set which is a directed
system, so that for every $a, b \in A$ there is a $c \in A$ such that
$a, b \preceq c$.  Also let $(X, \mathcal{A}, \mu)$ be a probability
space, and for each element $a$ of the directed set $A$, let
$\mathcal{B}_a$ be a $\sigma$-subalgebra of the measurable sets in
$X$.  Suppose that these $\sigma$-algebras are compatible with the
ordering on $A$, in the sense that
\begin{equation}
\label{mathcal{B}_a subseteq mathcal{B}_b}
        \mathcal{B}_a \subseteq \mathcal{B}_b
\end{equation}
when $a, b \in A$ satisfy $a \preceq b$.  This includes the usual
notion of a filtration on $X$, and some aspects of martingales still
work in this setting, as in \cite{ks}.  However, standard results
about maximal functions do not always hold, for instance.

        Let $I$ be an infinite set, and let $(X_j, \mathcal{A}_j, \mu_j)$
be a probability space for each $j \in I$.  Consider the product $X =
\prod_{j \in I} X_j$, with the usual product $\sigma$-algebra
$\mathcal{A}$ of measurable sets and probability measure $\mu$.  If
$K$ is a nonempty finite subset of $I$, then let $\mathcal{B}_K$ be
the $\sigma$-subalgebra of measurable subsets of $X$ that correspond
to a product of measurable subset of $\prod_{j \in K} X_j$ with
$\prod_{j \in I \setminus K} X_j$.  Equivalently, $\mathcal{B}_K$
consists of the inverse images of measurable subsets of $\prod_{j \in
  K} X_j$ under the obvious coordinate projection from $X$ onto
$\prod_{j \in K} X_j$.  If $L$ is another finite subset of $I$ that
contains $K$, then $\mathcal{B}_K \subseteq \mathcal{B}_L$.  The
collection of nonempty finite subsets of $I$ is a directed system with
respect to inclusion, and this defines a compatible family of
$\sigma$-subalgebras of measurable sets in $X$.  In this situation,
the usual arguments about maximal functions do not work, even when
$I$ is countable.  In particular, there are problems with pointwise
convergence.

        Now let $\preceq$ be the partial ordering on ${\bf Z}_+$ where
$a \preceq b$ when $b$ is an integer multiple of $a$.  Of course,
${\bf Z}_+$ is a directed system with respect to this ordering.
If $r = \{r_j\}_{j = 1}^\infty$ is a sequence of integers with $r_j
\ge 2$ for each $j$, and if $R_l$ is as in (\ref{R_l = prod_{j = 1}^l
  r_j, 2}), then the set of $R_l$'s is linearly-ordered with respect to
$\preceq$, and hence is a directed system.  Let $E$ be an infinite
subset of ${\bf Z}_+$, and suppose that $E$ is also a directed system
with respect to this ordering.  Let $p$ be a prime number, and let
$c_E(p)$ be the supremum of the nonnegative integers $k$ for which
there is an $R \in E$ that is an integer multiple of $p$.  Thus
$c_E(p)$ is either a nonnegative integer or $+\infty$ for each $p$.
Because $E$ is infinite, either $c_E(p) = +\infty$ for some $p$, or
$c_E(p) > 0$ for infinitely many $p$.  Using $E$, we get a
translation-invariant topology on ${\bf Z}$, for which a local base
for the topology at $0$ is given by the sets $R \, {\bf Z}$ with $R
\in E$.  As usual, ${\bf Z}$ is a topological ring with respect to
this topology.  If $r$ is as before and $c_r(p) = c_E(p)$ for every
prime number $p$, then the topology on ${\bf Z}$ corresponding to $E$
is the same as the $r$-adic topology discussed previously.

        Consider the Cartesian product
\begin{equation}
\label{X_E = prod_{R in E} ({bf Z} / R {bf Z})}
        X_E = \prod_{R \in E} ({\bf Z} / R \, {\bf Z}).
\end{equation}
This is a compact commutative ring with respect to coordinatewise
addition and multiplication, and using the product topology associated
to the discrete topology on ${\bf Z} / R \, {\bf Z}$ for each $R \in
E$.  There is a natural homomorphism from ${\bf Z}$ into $X_E$,
defined by the canonical quotient mappings from ${\bf Z}$ onto ${\bf
  Z} / R \, {\bf Z}$ for each $R \in E$.  It is easy to see that the
kernel of this homomorphism is trivial, because $E$ is infinite.  This
homomorphism is also a homeomorphism from ${\bf Z}$ onto its image in
$X_E$ with respect to the topology on ${\bf Z}$ associated to $E$ as
in the preceding paragraph.

        Let $x = \{x_R\}_{R \in E}$ be an element of $X_E$, so that $x_R
\in {\bf Z} / R \, {\bf Z}$ for each $R \in E$.  If $R_1, R_2 \in E$ and 
$R_1 \preceq R_2$, then $R_2 \, {\bf Z} \subseteq R_1 \, {\bf Z}$, and we 
get a natural ring homomorphism from ${\bf Z} / R_2 \, {\bf Z}$ onto
${\bf Z} / R_1 \, {\bf Z}$.  Let us say that $x \in X_E$ is coherent if
the $x_{R_1}$ is the image of $x_{R_2}$ under the natural homomorphism
from ${\bf Z} / R_2 \, {\bf Z}$ onto ${\bf Z} / R_1 \, {\bf Z}$ for
every $R_1, R_2 \in E$ such that $R_1 \preceq R_2$.  The set of coherent
elements of $X_E$ forms a closed sub-ring of $X_E$.  The natural homomorphism
from ${\bf Z}$ into $X_E$ maps ${\bf Z}$ into the set of coherent elements
of $X_E$, and in fact the set of coherent elements of $X_E$ is the
same as the closure of the image of ${\bf Z}$ in $X_E$.  Let ${\bf Z}_E$
be the set of coherent elements of $X_E$.  If $r = \{r_j\}_{j = 1}^\infty$
is as before and $c_r(p) = c_E(p)$ for every prime number $p$, then
${\bf Z}_E$ can be identified with ${\bf Z}_r$.

        There is a natural homomorphism from ${\bf Z}_E$ onto 
${\bf Z} / R \, {\bf Z}$ for each $R \in E$, which is the restriction
to ${\bf Z}_E$ of the coordinate mapping from $X_E$ onto 
${\bf Z} / R \, {\bf Z}$.  Let $\mathcal{B}_R$ be the collection of
subsets of ${\bf Z}_E$ that can be expressed as the inverse image of
a subset of ${\bf Z} / R \, {\bf Z}$ under the mapping just defined,
for each $R \in E$.  This is a $\sigma$-subalgebra of the Borel sets 
in ${\bf Z}_E$.  If $R_1, R_2 \in E$ satisfy $R_1 \preceq R_2$,
then it is easy to see that $\mathcal{B}_{R_1} \subseteq \mathcal{B}_{R_2}$.
Thus we get a family of $\sigma$-subalgebras of the Borel sets in ${\bf Z}_E$
indexed by $E$ and compatible with the ordering on $E$.

        Similarly,
\begin{equation}
\label{Y_E = prod_{R in E} ({bf R} / R {bf Z})}
        Y_E = \prod_{R \in E} ({\bf R} / R \, {\bf Z})
\end{equation}
is a compact commutative topological group with respect to
coordinatewise addition and the product topology associated to the
usual quotient topology on ${\bf R} / R \, {\bf Z}$ for each $R \in
E$.  There is a natural continuous group homomorphism from ${\bf R}$
as a commutative topological group with respect to addition into
$Y_E$, defined by the canonical quotient mapping from ${\bf R}$ onto
${\bf R} / R \, {\bf Z}$ for each $R \in E$.  If $R_1, R_2 \in E$ and
$R_1 \prec R_2$, then $R_2 \, {\bf Z} \subseteq R_1 \, {\bf Z}$, which
leads to a continuous group homomorphism from ${\bf R} / R_2 \, {\bf
  Z}$ onto ${\bf R} / R_1 \, {\bf Z}$.  As usual, an element $y =
\{y_R\}_{R \in E}$ of $Y$ is said to be coherent if $y_{R_1}$ is the
image of $y_{R_2}$ under the natural mapping from ${\bf R} / R_2 \,
{\bf Z}$ onto ${\bf R} / R_1 \, {\bf Z}$ for every $R_1, R_2 \in E$
with $R_1 \preceq R_2$.  The set of coherent elements of $Y_E$ is a
closed subgroup of $Y_E$, which can be identified with the quotient of
${\bf R} \times {\bf Z}_E$ by the image of ${\bf Z}$ under the obvious
diagonal embedding.  Let $C_R$ be the collection of subsets of the
group of coherent elements of $Y_E$ that can be expressed as the
inverse image of a Borel set in ${\bf R} / R \, {\bf Z}$ under the
corresponding coordinate mapping, for each $R \in E$.  If $R_1, R_2
\in E$ and $R_1 \prec R_2$, then $C_{R_1} \subseteq C_{R_2}$, as
before.  This defines a family of $\sigma$-subalgebras of the Borel
sets in the group of coherent elements of $Y_E$ indexed by $E$
which is compatible with the ordering on $E$.


\begin{thebibliography}{49}

\addcontentsline{toc}{section}{References}



\bibitem {a-m} M.~Atiyah and I.~MacDonald, {\it Introduction to
Commutative Algebra}, Addison-Wesley, 1969.

\bibitem {b-j-l-p-s} S.~Bates, W.~Johnson, J.~Lindenstrauss,
D.~Preiss, and G.~Schechtman, {\it Affine approximation of Lipschitz
functions and nonlinear quotients}, Geometric and Functional Analysis
{\bf 9} (1999), 1092--1127.

\bibitem {b-l} Y.~Benyamini and J.~Lindenstrauss, {\it Geometric
Nonlinear Functional Analysis}, American Mathematical Society, 2000.

\bibitem {b-mac} G.~Birkhoff and S.~Mac Lane, {\it A Survey of Modern
  Algebra}, 4th edition, Macmillan, 1977.

\bibitem {cas} J.~Cassels, {\it Local Fields}, Cambridge University
Press, 1986.

\bibitem {c-w-1} R.~Coifman and G.~Weiss, {\it Analyse Harmonique
Non-Commutative sur Certains Espaces Homog\`enes}, Lecture Notes in
Mathematics {\bf 242}, Springer-Verlag, 1971.

\bibitem {c-w-2} R.~Coifman and G.~Weiss, {\it Extensions of Hardy
spaces and their use in analysis}, Bulletin of the American
Mathematical Society {\bf 83} (1977), 569--645.

\bibitem {d-s-1} G.~David and S.~Semmes, {\it Fractured Fractals and
Broken Dreams: Self-Similar Geometry through Metric and Measure},
Oxford University Press, 1997.

\bibitem {d-s-2} G.~David and S.~Semmes, {\it Regular mappings between
dimensions}, Publicacions Matem\`atiques {\bf 44} (2000), 369--417.

\bibitem {fa1} K.~Falconer, {\it The Geometry of Fractal Sets},
Cambridge University Press, 1986.

\bibitem {fa2} K.~Falconer, {\it Fractal Geometry: Mathematical
Foundations and Applications}, 2nd edition, Wiley, 2003.

\bibitem {f1} G.~Folland, {\it A Course in Abstract Harmonic
Analysis}, CRC Press, 1995.

\bibitem {f2} G.~Folland, {\it Real Analysis: Modern Techniques and
their Applications}, 2nd edition, Wiley, 1999.

\bibitem {gh} F.~Gehring, {\it The $L^p$-integrability of the partial
derivatives of a quasiconformal mapping}, Acta Mathematica {\bf 130}
(1973), 265--277.

\bibitem {gv} F.~Gouv\^ea, {\it $p$-Adic Numbers: An Introduction},
2nd edition, Springer-Verlag, 1997.

\bibitem {h1} J.~Heinonen, {\it Lectures on Analysis on Metric
Spaces}, Springer-Verlag, 2001.

\bibitem {h2} J.~Heinonen, {\it Geometric embeddings of metric
spaces}, Reports of the Department of Mathematics and Statistics {\bf
90}, University of Jyv\"askyl\"a, 2003.

\bibitem {h-r} E.~Hewitt and K.~Ross, {\it Abstract Harmonic
Analysis}, Volumes I, II, Springer-Verlag, 1970, 1979.

\bibitem {h-s} E.~Hewitt and K.~Stromberg, {\it Real and Abstract
Analysis}, Springer-Verlag, 1975.

\bibitem {h-w} W.~Hurewicz and H.~Wallman, {\it Dimension Theory},
Princeton University Press, 1969.

\bibitem {kat} Y.~Katznelson, {\it An Introduction to Harmonic
Analysis}, 3rd edition, Cambridge University Press, 2004.

\bibitem {ak1} A.~Knapp, {\it Basic Real Analysis}, Birkh\"auser, 2005.

\bibitem {ak2} A.~Knapp, {\it Advanced Real Analysis}, Birkh\"auser, 2005.

\bibitem {sk} S.~Krantz, {\it A Panorama of Harmonic Analysis},
Mathematical Association of America, 1999.

\bibitem {mac-b} S.~Mac Lane and G.~Birkhoff, {\it Algebra}, 3rd
  edition, AMS Chelsea, 1999.

\bibitem {mat} P.~Mattila, {\it Geometry of Sets and Measures in
Euclidean Spaces: Fractals and Rectifiability}, Cambridge University
Press, 1995.

\bibitem {mc} G.~McCarty, {\it Topology: An Introduction with
Application to Topological Groups}, 2nd edition, Dover, 1988.

\bibitem {r1} W.~Rudin, {\it Principles of Mathematical Analysis}, 3rd
edition, McGraw-Hill, 1976.

\bibitem {r2} W.~Rudin, {\it Real and Complex Analysis}, 3rd edition,
McGraw-Hill, 1987.

\bibitem {r3} W.~Rudin, {\it Fourier Analysis on Groups}, Wiley, 1990.

\bibitem {r4} W.~Rudin, {\it Functional Analysis}, 2nd edition,
McGraw-Hill, 1991.

\bibitem {s1} S.~Semmes, {\it An introduction to analysis on metric
spaces}, Notices of the American Mathematical Society {\bf 50} (2003),
438--443.

\bibitem {s2} S.~Semmes, {\it Some Novel Types of Fractal Geometry},
Oxford University Press, 2001.

\bibitem {st1} E.~Stein, {\it Singular Integrals and Differentiability
Properties of Functions}, Princeton University Press, 1970.

\bibitem {st2} E.~Stein, {\it Topics in Harmonic Analysis Related to
  the Littlewood--Paley Theory}, Annals of Mathematics Studies {\bf
  63}, Princeton University Press, 1970.

\bibitem {st3} E.~Stein, {\it Harmonic Analysis: Real-Variable Methods,
Orthogonality, and Oscillatory Integrals}, with the assistance of
T.~Murphy, Princeton University Press, 1993.

\bibitem {st-sh-1} E.~Stein and R.~Shakarchi, {\it Fourier Analysis:
An Introduction}, Princeton University Press, 2003.

\bibitem {st-sh-2} E.~Stein and R.~Shakarchi, {\it Real Analysis:
Measure Theory, Integration, and Hilbert Spaces}, Princeton University
Press, 2005.

\bibitem {st-sh-3} E.~Stein and R.~Shakarchi, {\it Functional
  Analysis}, Princeton University Press, 2011.

\bibitem {s-w} E.~Stein and G.~Weiss, {\it Introduction to Fourier
Analysis on Euclidean Spaces}, Princeton University Press, 1971.

\bibitem {ks} K.~Stromberg, {\it Probability for Analysts}, lecture
  notes prepared by K.~Ravindran, Chapman \& Hall, 1994.

\bibitem {ds} D.~Sullivan, {\it Linking the universalities of
Milnor--Thurston, Feigenbaum and Ahlfors--Bers}, in {\it Topological
Methods in Modern Mathematics}, 543--564, Publish or Perish, 1993.

\bibitem {t} M.~Taibleson, {\it Fourier Analysis on Local Fields},
Princeton University Press, 1975.

\bibitem {t-v} P.~Tukia and J.~V\"ais\"al\"a, {\it Quasisymmetric
embeddings of metric spaces}, Annales Academiae Scientiarum Fennicae
Series A I Mathematica {\bf 5} (1980), 97--114.

\bibitem {z} A.~Zygmund, {\it Trigonometric Series}, Volumes I, II,
3rd edition, with a foreword by R.~Fefferman.




\end{thebibliography}
\end{document}